\documentclass[11pt]{amsart}
\usepackage[colorlinks=true, pdfstartview=FitV, linkcolor=blue, citecolor=blue, urlcolor=blue, breaklinks=true]{hyperref}
\usepackage{amsmath,amsfonts,amssymb,amsthm,amscd,comment,euscript}
\usepackage{ytableau}
\usepackage{xcolor}
\usepackage{paralist}
\usepackage[vcentermath]{youngtab}
\usepackage[latin1]{inputenc}
\usepackage{tikz}
\usepackage[all]{xy}
\usepackage{tikz}
\usepackage{array}
\newtheorem{cor}{Corollary}[subsection]
\newtheorem{lem}{Lemma}[subsection]
\newtheorem{prop}{Proposition}[subsection]

\newcommand{\nc}{\newcommand}

\newcounter{cnt}
 \makeatletter
\def\mydggeometry{\makeatletter\dg@YGRID=1\dg@XGRID=20\unitlength=0.003pt\makeatother}
\makeatother \theoremstyle{remark}
\numberwithin{equation}{section}

\let\bwdg\bigwedge
\def\bigwedge{{\textstyle\bwdg}}
\theoremstyle{definition}
\newtheorem*{defn}{Definition}
\newtheorem*{ex}{Example}
\theoremstyle{definition}
\newtheorem{thm}{Theorem}[subsection]

\newtheorem{rem}{Remark}[subsection]

\newcommand\al{\alpha}
\newcommand{\ts}{\textstyle}
\newcommand\Lg{\mathfrak{g}}
\newcommand\Lh{\mathfrak{h}}
\newcommand\Lb{\mathfrak{b}}
\newcommand\Ln{\mathfrak{n}}
\newcommand\U{\mathfrak{U}}

\DeclareMathOperator{\height}{ht}
\nc{\C}{\mathbb C }
\nc{\D}{\mathbb D }
\nc{\Z}{\mathbb Z }
\nc{\N}{\mathbb N }
\nc{\R}{\mathbb R }
\nc{\Q}{\mathbb Q }
\newcommand{\wt}{\operatorname{wt}}

\begin{document}

\title[Realization of affine type A KR-crystals via polytopes]{Realization of affine type A Kirillov-Reshetikhin crystals via polytopes}
\author{Deniz Kus}
\address{Deniz Kus:\newline
Mathematisches Institut, Universit\"at zu K\"oln, Germany}
\email{dkus@math.uni-koeln.de}
\thanks{The author was sponsored by the "SFB/TR 12 - Symmetries and Universality in Mesoscopic Systems".}

\begin{abstract}
On the polytope defined in \cite{FFL11}, associated to any rectangle highest weight, we define a structure of an type $A_n$-crystal. We show, by using the Stembridge axioms, that this crystal is isomorphic to the one obtained from Kashiwara's crystal bases theory. Further we define on this polytope a bijective map and show that this map satisfies the properties of a weak promotion operator. This implies in particular that we provide an explicit realization of Kirillov-Reshetikhin crystals for the affine type $A^{(1)}_n$ via polytopes.
\end{abstract}
\maketitle \thispagestyle{empty}
%%%%%%%%%%%%%%%%%%%%%%%%%%%%%%%%%%%%%%%%%%%%%%%%%%%%%%%%%%%%%%%%%%%%%%%%%%%%%%%%%%%%%%%%%%%%%%%%%%%%%%%%%%%%%%%%%%%%%%%%%%%%%%%%%%%
%         Introduction
%%%%%%%%%%%%%%%%%%%%%%%%%%%%%%%%%%%%%%%%%%%%%%%%%%%%%%%%%%%%%%%%%%%%%%%%%%%%%%%%%%%%%%%%%%%%%%%%%%%%%%%%%%%%%%%%%%%%%%%%%%%%%%%%%%%
\section{Introduction}
Let $\Lg$ be a affine Lie algebra and $\U^{'}_q(\Lg)$ be the corresponding quantum algebra without derivation. The irreducible representations are classified in \cite{CP95},\cite{CP98} in terms of Drinfeld polynomials. A certain subclass of these modules, that gained a lot of attraction during the last decades, are the so called Kirillov-Reshetikhin modules $KR(m,\omega_i,a)$, where $i$ is a node in the classical Dynkin diagram and $m$ is a positive integer. One of the main tools for studying such representations is Kashiwara's crystal bases theory \cite{Kashi91}. This theory was originally defined for representations for $\U_q(\Lg)$, however it can be nonetheless defined in the setting of $\U^{'}_q(\Lg)$ modules, respecting that crystal bases might not always exist. It was first conjectured in \cite{HKOTT02}, that $KR(m,\omega_i,a)$ admits a crystal bases and this 
was proven in type $A^{(1)}_n$ in \cite{KKMMNN92} and in all non-exceptional cases in \cite{O07},\cite{OS08}. We denote this crystal by $KR^{m,i}$ and call it a Kirillov-Reshetikhin crystal.\par A promotion operator $pr$ on a crystal $B$ of type $A_n$ is defined to be a map satisfying several conditions, namely that $pr$ shifts the content, $pr\circ\tilde{e}_j=\tilde{e}_{j+1}\circ pr$, $pr\circ\tilde{f}_j=\tilde{f}_{j+1}\circ pr$ for all $j\in\{1,\cdots,n-1\}$ and $pr^{n+1}=id$, where $\tilde{e}_j$ and $\tilde{f}_j$ respectively are the Kashiwara operators. If the latter condition is not satisfied, but $pr$ is still bijective, then the map $pr$ is called a weak promotion operator (see also \cite{BST10}). The advantage of such (weak) promotion oparators are that we can associate to a given crystal $B$ of type $A_n$ a (weak) affine crystal by setting $\tilde{f}_0:=pr^{-1}\circ \tilde{f}_1\circ pr, \mbox{ and } \tilde{e}_0:=pr^{-1}\circ \tilde{e}_1\circ pr$.\par On the set of all semi-standard Young tableaux of rectangle shape, which is a realization of $B(m\omega_i)$ the $\U_q(A_n)$-crystal associated to the irreducible module of highest weight $m\omega_i$, Schützenberger defined a promotion operator $pr$, called the Schützenberger's promotion operator \cite{S72}, which is the analogue of the cyclic Dynkin diagram automorphism $i\mapsto i+1\mod (n+1)$ on the level of crystals, by using  jeu-de-taquin. Given a tableaux $T$ over the alphabet $1\prec 2\cdots \prec n+1$, $pr(T)$ is obtained from $T$ by removing all letters $n+1$, adding one to each entry in the remaining tableaux, using jeu-de-taquin to slide all letters up and finally filling the holes with 1's (see also Section~\ref{section6}). One of the combinatorial descriptions of $KR^{m,i}$ in the affine $A^{(1)}_n$ type was provided by Schimozono in \cite{S02}. It was shown that, as a $\{1,\cdots,n\}$-crystal, $KR^{m,i}$ is isomorphic to $B(m\omega_i)$ and the affine crystal constructed from $B(m\omega_i)$ using Schützenberger's promotion operator is isomorphic to the Kirillov-Reshetikhin crystal $KR^{m,i}$. The two ways of computing the affine crystal structure, one given by \cite{KKMMNN92} and the other by \cite{S02}, are shown to be equivalent in \cite{OSS03}. Another combinatorial model in this type without using a promotion operator is described in \cite{KW12}. In this paper, we introduce a new realization of Kirillov-Reshetikhin crystals of type $A^{(1)}_n$.\par In \cite{FFL11} the authors have constructed for all dominant integral $A_n$ weights $\lambda$ a polytope in $\R^{n\frac{n-1}{2}}$ and a basis of the irreducible $A_n$ module of highest weight $\lambda$ and have shown that the basis elements are parametrized by the integral points. For $\lambda=m\omega_i$ we can understand this polytope in $\R^{i(n-i+1)}$ and denote the intersection of this polytope with $\Z_{+}^{i(n-i+1)}$ by $B^{m,i}$. We define certain maps on $B^{m,i}$ and show that this becomes a crystal of type $A_n$. As a set, we can identify $B^{m,i}$ with certain blocks of height $n-i+1$ and width $i$
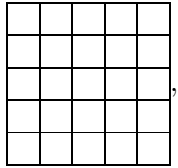
$$\yng(5,5,5,5,5),$$ where the boxes are filled, under some assumptions, with some non-negative integers (see Definition~\ref{maindef}). The crystal $B^{m,i}$ has no known explicit combinatorial bijection to other combinatorial models of  crystals induced by representations, such as the Young tableaux model \cite{NK94} or the set of certain Nakajima monomials \cite{Nak03}, which makes an isomorphism very complicated. Using the realization of crystal bases via Nakajima monomials, we can construct certain local $A_2$ isomorphisms on our underlying polytope $B^{m,i}$ and prove that the so called Stembridge axioms are satisfied. These axioms precisely characterize the set of crystals of representations in the class of all crystals. Our first important theorem is therefore the following:
\vskip 5pt\noindent
{\bf Theorem~A}.\
%\ref{basis}}
{\it The polytope $B^{m,i}$ is as an $A_n$ crystal isomorphic to $B(m\omega_i).$}
\vskip 5pt
In order to obtain an (weak) affine crystal structure we define a map $pr$ on $B^{m,i}$, which is given by an algorithm consisting of $i$ steps (see (\ref{algorithm})), and show that this map satisfies the conditions for a weak promotion operator. In particular, this implies that this map $pr$ is the unique promotion operator on $B^{m,i}\cong B(m\omega_i)$ and the polytope becomes an affine crystal. To be more precise, we prove the following main theorem of our paper:
\vskip 5pt\noindent
{\bf Theorem~B}.\
%\ref{basis}}
{\it The associated affine crystal $B^{m,i}$ using $pr$ is isomorphic to the Kirillov-Reshetikhin crystal $KR^{m,i}$.}
\vskip 5pt 
Our paper is organized as follows: in Section~\ref{section2} we fix some notation and present the main definitions, especially the definition of our polytope. In Section~\ref{section3} we equip our main object with a crystal structure. In  Section~\ref{section4} Nakajima monomials are recalled and in Section~\ref{section5} Theorem A is proven. Finally, in Section~\ref{section6} the promotion operator is defined by an algorithm and the corresponding affine crystal is identified with the KR crystal, proving Theorem B.
\vskip12pt
\textbf{Acknowledgements:} 
The author would like to thank Ghislain Fourier and Peter Littelmann for their helpful discussions and Vyjayanthi Chari and the University of California at Riverside for their hospitality during his stays there, when some of the ideas of the current paper were developed.
%%%%%%%%%%%%%%%%%%%%%%%%%%%%%%%%%%%%%%%%%%%%%%%%%%%%%%%%%%%%%%%%%%%%%%%%%%%%%%%%%%%%%%%%%%%%%%%%%%%%%%%%%%%%%%%%%%%%%%%%%%%%%%%%%%%
%         
%%%%%%%%%%%%%%%%%%%%%%%%%%%%%%%%%%%%%%%%%%%%%%%%%%%%%%%%%%%%%%%%%%%%%%%%%%%%%%%%%%%%%%%%%%%%%%%%%%%%%%%%%%%%%%%%%%%%%%%%%%%%%%%%%%%
\section{Notation and main definitions}\label{section2}
Let $\Lg$ be a complex affine Lie algebra of rank $n$ and fix a Cartan subalgebra $\Lh$ in $\Lg$ and a Borelsubalgebra $\Lb\supseteq\Lh$. We denote by $\Phi\subseteq \Lh^*$ the root system of the Lie algebra, and, corresponding to the choice of $\Lb$ let $\Phi^+$ be the subset of positive roots. Further, we denote by $\Pi=\{\alpha_0,\cdots,\alpha_n\}$ the corresponding basis of $\Phi$ and the basis of the dual root system $\Phi^{\vee}\subseteq \Lh$ is denoted by $\Pi^\vee=\{\alpha_0^\vee,\cdots,\alpha_n^\vee\}$. Let  $\Lg=\Ln^+\oplus\Lh\oplus \Ln^-$ be a Cartan decomposition and for a given root $\alpha\in \Phi$ let $\Lg_{\alpha}$ be the corresponding root space.
For a dominant integral weight $\lambda$ we denote by $V(\lambda)$ the irreducible $\Lg$-module with
highest weight $\lambda$. Fix a highest weight vector $v_\lambda\in V(\lambda)$,
then $V(\lambda)=\U(\Ln^-)v_\lambda$, where $\U(\Ln^-)$ denotes the universal
enveloping algebra of $\Ln^-$. For an indetermined element $q$ we denote by $\U^{'}_q(\Lg)$ be the corresponding quantum algebra without derivation. The irreducible representations are classified in \cite{CP95},\cite{CP98} in terms of Drinfeld polynomials. One of the major goals in representation theory is to find nice expressions for the character of objects in the category $\mathcal{O}^q_{int}$ (see \cite{HK02}). From the theory of crystal bases, introduced by Kashiwara in \cite{Kashi91}, we can compute the character of a given module $M$ as follows:
$$ ch M=\sum_{\mu}\sharp (B_{\mu}) e^{\mu},$$ whereby $(L,B)$ is the crystal bases of $M$ (see also \cite{HK02}). From now on let $\Lg$ be the affine Lie algebra $$A^{(1)}_n=\mathfrak{sl}_{n+1}\otimes \C[t,t^{-1}]\oplus \C c\oplus \C d,$$ with index set $\hat{I}=\{0,1,\cdots,n\}\supseteq I=\{1,\cdots,n\}$. Note that the classical positive roots are all of the form $$\alpha_{i,j}=\alpha_i+\alpha_{i+1}+\cdots+\alpha_j,\ \mbox{for $1\leq i\leq j\leq n$}.$$ Further let $P=\bigoplus_{i\in I}\Z\omega_i$ be the set of classical integral and $P^+=\bigoplus_{i\in I}\Z_{+}\omega_i$ be the set of classical dominant integral weights. In order to realize the crystal graph of the so called Kirillov-Reshetikhin modules $KR(m,\omega_i,a)$, for $i\in I, m\in \Z_{+}$, we will define now the underlying combinatorial model in this paper, which we will denote by $B^{m,i}$. For more details regarding KR-modules we refer to a series of papers (\cite{C01},\cite{CM06},\cite{FoL07}).

\subsection{The polytope \texorpdfstring{$B^{m,i}$}{B}}
In this subsection we will define the set $B^{m,i}$, our main object in this paper and discuss its combinatorics which is crucial for the realization of KR-crystals. 
\begin{defn}\label{maindef}
Let $B^{m,i}$ be the set of all following patterns:

$$\ytableausetup{boxsize=2.8em}
\begin{ytableau}
a_{1,i} & a_{2,i} & \dots
& a_{i-1,i}& a_{i,i} \\
a_{1,i+1} & a_{2,i+1} & \dots
& \scriptstyle a_{i-1,i+1}& a_{i,i+1} \\ 
\vdots & \vdots & \vdots & \vdots & \vdots \\
a_{1,n} & a_{2,n} & \dots
&  a_{i-1,n}& a_{i,n} \\\end{ytableau}$$

filled with non-negative integers, such that $\sum^k_{s=1}a_{\beta(s)}\leq m$ for all sequences \[
(\beta(1),\dots, \beta(k)), \ k\ge 1
\] satisfying the following: $\beta(1)=\al_{1,i}, \beta(k)=\al_{i,n}$ and if $\beta(s)=\al_{p,q}$ then the next element in the sequence is either of the form 
$\beta(s+1)=\al_{p,q+1}$  or $\beta(s+1)=\al_{p+1,q}.$
\end{defn}
\begin{ex}
$$
{\Yvcentermath1
\young(10,21,01) \mbox{ $\in B^{5,2}$, but } {\Yvcentermath1
\young(100,013,101)\mbox{ $\notin B^{5,3}$.}}}
$$
\end{ex}
\begin{rem}\mbox{}
\begin{enumerate}
\item For any element in $B^{m,i}$ the columns are numbered from $1$ to $i$ and the rows are numbered from $i$ to $n$.
\item A sequence $$
\mathbf{b}=(\beta(1),\dots, \beta(k)), \ k\ge 1
$$ satisfying the rule from Definition~\ref{maindef} is called a Dyck path. The notion of a Dyck path occurs already in \cite{FFL11}. We will denote the set of all such paths by $\mathbf D$.
\end{enumerate} 
\end{rem}
\begin{rem}
Note that the name ``polytope" is justified, since $B^{m,i}$ reflects the integral points of some polytope in $\R^{i(n-i+1)}$:$$B^{m,i}\cong \{(a_{r,s})\in\R^{i(n-i+1)}|\sum^k_{s=1}a_{\beta(s)}\leq m, \mbox{ for all $\mathbf b\in\mathbf{D}$ }\}\cap \Z_+^{i(n-i+1)}.$$
\end{rem}
What we want to show now is that the set $B^{m,i}$ carries an (affine) crystal structure. Moreover, our goal is to show that this is exactly the crystal graph of the KR-module $KR(m,\omega_i,a)$, i.e. we have an isomorphism of crystals. The strongest indication that this conjecture might be true is the following modified result due to \cite{FFL11}.

\begin{thm}
$$\dim V(m\omega_i)=\sharp B^{m,i}$$
\end{thm}
%%%%%%%%%%%%%%%%%%%%%%%%%%%%%%%%%%%%%%%%%%%%%%%%%%%%%%%%%%%%%%%%%%%%%%%%%%%%%%%%%%%%%%%%%%%%%%%%%%%%%%%%%%%%%%%%%%%%%%%%%%%%%%%%%%%
%         
%%%%%%%%%%%%%%%%%%%%%%%%%%%%%%%%%%%%%%%%%%%%%%%%%%%%%%%%%%%%%%%%%%%%%%%%%%%%%%%%%%%%%%%%%%%%%%%%%%%%%%%%%%%%%%%%%%%%%%%%%%%%%%%%%%%
\section{Crystal Structure on \texorpdfstring{$B^{m,i}$}{B}}\label{section3}
With the purpose to show that we have a crystal structure on $B^{m,i}$, which is induced from a module we will first chop all necessary conditions of an abstract crystal. Let us start by giving the definition:

\subsection{Abstract crystals}
\begin{defn}\label{abcrystal}
Let $\hat{I}$ be a finite index set and let $A=(a_{i,j})_{i,j\in \hat{I}}$ be a generalized Cartan matrix with the Cartan datum $(A,\Pi,\Pi^{\vee},P,P^{\vee})$. A crystal associated with the Cartan datum $(A,\Pi,\Pi^{\vee},P,P^{\vee})$ is a set $B$ together with the maps $\wt:B \rightarrow P$, $\tilde{e}_l,\tilde{f}_l: B\rightarrow B \cup \{0\}$, and $\epsilon_l,\varphi_l:B\rightarrow \Z\cup \{-\infty\}$ satisfying the following properties for all $l\in \hat{I}$:
\begin{enumerate}
\item $\varphi_l(b)=\epsilon_l(b)+\langle \alpha_l^\vee,\wt(b)\rangle$
\item $\wt(\tilde{e}_{l}b)=\wt(b)+\alpha_l$ if $\tilde{e}_{l}b\in B$
\item $\wt(\tilde{f}_{l}b)=\wt(b)-\alpha_l$ if $\tilde{f}_{l}b\in B$
\item $\epsilon_l(\tilde{e}_{l}b)=\epsilon_l(b)-1$, $\varphi_l(\tilde{e}_{l}b)=\varphi_l(b)+1$ if $\tilde{e}_{l}b\in B$
\item $\epsilon_l(\tilde{f}_{l}b)=\epsilon_l(b)+1$, $\varphi_l(\tilde{f}_{l}b)=\varphi_l(b)-1$ if $\tilde{f}_{l}b\in B$
\item $\tilde{f}_{l}b=b'$ if and only if $\tilde{e}_{l}b'=b$ for $b,b'\in B$
\item if $\varphi_l(b)=-\infty$ for $b\in B$, then $\tilde{f}_{l}b=\tilde{e}_{l}b$=0.
\end{enumerate}
Further a crystal $B$ is said to be semiregular if the equalities $$\epsilon_l(b)=\max\{k\geq0|\tilde{e}_{l}^{k}b\neq 0\},\quad \varphi_l(b)=\max\{k\geq0|\tilde{f}_{l}^{k}b\neq 0\}$$ hold.
\end{defn}

Hence our aim is to define the Kashiwara operators $\tilde{f}_{l}$ and $\tilde{e}_{l}$, which will act on the set $B^{m,i}$ for all $l\in\{0,\cdots,n\}$, such that the properties in Definition~\ref{abcrystal} are fulfilled. Furthermore we will show that $B^{m,i}$ becomes a semiregular crystal. Our strategy is as follows: first we are going to define a classical crystal structure on $B^{m,i}$, which in particular means that we define the Kashiwara operators for all $l\in\{1,\cdots,n\}$. Subsequently, we show that this is precisely the crystal graph of the irreducible module $V(m\omega_i)$ and then we exploit the existence of Schützenberger's promotion operator to define the Kashiwara operators for the node $0$. For more details regarding the Schützenberger promotion operator we refer to \cite{S72} or Section~\ref{section6}.

\subsection{Crystal Structure on \texorpdfstring{$B^{m,i}$}{B}}

As already mentioned, our aim in this section is to define the maps $\wt,\tilde{e}_{l},\tilde{f}_{l},\epsilon_l,\varphi_l$ for all $l\in I$ as in the Definition~\ref{abcrystal}, such that the properties (1)-(7) are fulfilled. So let $A$ be an arbitrary element in $B^{m,i}$, then we define \begin{equation}\label{weight}\wt(A)=m\omega_i-\sum_{1\leq p\leq i, i\leq q\leq n}a_{p,q}\alpha_{p,q}.\end{equation}
In order to define what the Kashiwara operators are, we need much more spadework. In the following we define some useful maps and integers, such that these integers will completely determine the rule at which ``place" the action is given. So we define the maps $\varphi_l,\epsilon_l: B^{m,i}\longrightarrow \Z_{\geq 0}$ for all $l\in I$ by:
\begin{equation}\label{ttt1}\varphi_l(A)=\begin{cases}
  m-\sum^{i-1}_{j=1} a_{j,i}-\sum^{n}_{j=i} a_{i,j},  & \text{if $l=i$}\\
  \sum^{p^l_+(A)}_{j=1}a_{j,l-1}-\sum^{p^l_+(A)-1}_{j=1}a_{j,l}, & \text{if $l>i$}\\
  \sum^{n}_{j=p^l_{-}(A)}a_{l+1,j}-\sum^{n}_{j=p^l_{-}(A)+1}a_{l,j}, & \text{if $l<i$}
\end{cases}\end{equation}
\begin{equation}\label{ttt2}\epsilon_l(A)=\begin{cases}
  a_{i,i},  & \text{if $l=i$}\\
  \sum^{i}_{j=q^l_+(A)}a_{j,l}-\sum^{i}_{j=q^l_{+}(A)+1}a_{j,l-1}, & \text{if $l>i$}\\
  \sum^{q^l_{-}(A)}_{j=i}a_{l,j}-\sum^{q^l_{-}(A)-1}_{j=i}a_{l+1,j}, & \text{if $l<i$,}
\end{cases}\end{equation}
whereby 
\begin{equation}\label{1}p^l_+(A)=\min\{1\leq p\leq i|\sum^p_{j=1}a_{j,l-1}+\sum^{i}_{j=p}a_{j,l}=\max_{1\leq q\leq i}\{\sum^q_{j=1}a_{j,l-1}+\sum^{i}_{j=q}a_{j,l}\}\}\end{equation}
\begin{equation}\label{2}q^l_+(A)=\max\{1\leq p\leq i|\sum^p_{j=1}a_{j,l-1}+\sum^{i}_{j=p}a_{j,l}=\max_{1\leq q\leq i}\{\sum^q_{j=1}a_{j,l-1}+\sum^{i}_{j=q}a_{j,l}\}\}\end{equation}
\begin{equation}\label{3}p^l_-(A)=\max\{i\leq p\leq n|\sum^p_{j=i}a_{l,j}+\sum^{n}_{j=p}a_{l+1,j}=\max_{i\leq q\leq n}\{\sum^q_{j=i}a_{l,j}+\sum^{n}_{j=q}a_{l+1,j}\}\}\end{equation}
\begin{equation}\label{4}q^l_-(A)=\min\{i\leq p\leq n|\sum^p_{j=i}a_{l,j}+\sum^{n}_{j=p}a_{l+1,j}=\max_{i\leq q\leq n}\{\sum^q_{j=i}a_{l,j}+\sum^{n}_{j=q}a_{l+1,j}\}\}.\end{equation}
\begin{rem}\label{notrem}
Note that the integers (\ref{ttt1})-(\ref{ttt2}) for $l\neq i$ and (\ref{1})-(\ref{4}) depend only on two given columns or two given rows of $A$. Therefore one can define these integers for any given two columns or rows $\mathbf{a}$ and $\mathbf{b}$ and denote them alternatively by $p_{\pm}(\mathbf{a},\mathbf{b})$,$q_{\pm}(\mathbf{a},\mathbf{b})$ and $\epsilon(\mathbf{a},\mathbf{b})$, $\varphi(\mathbf{a},\mathbf{b})$ respectively. For instance we will use in some places the notation $q_{-}(\mathbf{a_l},\mathbf{a_{l+1}})$ instead of (\ref{4}), if $\mathbf{a_l}$ and $\mathbf{a_{l+1}}$ is the $l$-th column and $(l+1)$-th column respectively of $A$, and $\epsilon(\mathbf{a_l},\mathbf{a_{l+1}})$ instead of (\ref{ttt2}).
\end{rem}
The first fact we want to note about these maps is the following lemma:
\begin{lem}\label{fest}
The map $\varphi_l$ is uniquely determined by the map $\epsilon_l$ and conversely the map $\epsilon_l$ is uniquely determined by the map $\varphi_l$. Particularly we have $$\varphi_l(A)=\epsilon_l(A)+\langle \alpha_l^\vee,\wt(A)\rangle.$$
\proof
Assume $A$ is an arbitrary element and let $p^l_+(A)$, $p^l_-(A)$, $q^l_+(A)$, $q^l_-(A)$ be the integers described in (\ref{1})-(\ref{4}). Since the statement is obvious for $l=i$ we presume $l\neq i$. Then, because of  $$\sum^{p^l_{+}(A)}_{j=1}a_{j,l-1}+\sum^{i}_{j=p^l_+(A)}a_{j,l}=\sum^{q^l_+(A)}_{j=1}a_{j,l-1}+\sum^{i}_{j=q^l_+(A)}a_{j,l}$$
and
$$\sum^{p^l_{-}(A)}_{j=i}a_{l,j}+\sum^{n}_{j=p^l_-(A)}a_{l+1,j}=\sum^{q^l_{-}(A)}_{j=i}a_{l,j}+\sum^{n}_{j=q^l_-(A)}a_{l+1,j},$$
it follows that \begin{equation}\label{was}\sum^{q^l_{+}(A)}_{j=p^l_{+}(A)+1}a_{j,l-1}-\sum^{q^l_{+}(A)-1}_{j=p^l_{+}(A)}a_{j,l}=\sum^{p^l_{-}(A)}_{j=q^l_{-}(A)+1}a_{l,j}-\sum^{p^l_{-}(A)-1}_{j=q^l_{-}(A)}a_{l+1,j}=0.\end{equation} Therefore, if $l>i$ we arrive at \begin{flalign*}\epsilon_{l}(A)&=\sum^{i}_{j=q^l_+(A)}a_{j,l}-\sum^{i}_{j=q^l_{+}(A)+1}a_{j,l-1}=\sum^{i}_{j=p^l_+(A)}a_{j,l}-\sum^{i}_{j=p^l_{+}(A)+1}a_{j,l-1}&\\&=\sum^{p^l_+(A)}_{j=1}a_{j,l-1}-\sum^{p^l_+(A)-1}_{j=1}a_{j,l}-(\sum^{i}_{j=1}a_{j,l-1}-\sum^{i}_{j=1}a_{j,l})=\varphi_l(A)-\langle \alpha_l^\vee,\wt(A)\rangle\end{flalign*} and if $l<i$ we obtain again with (\ref{was}) \begin{flalign*}\epsilon_{l}(A)&=\sum^{p^l_{-}(A)}_{j=i}a_{l,j}-\sum^{p^l_{-}(A)-1}_{j=i}a_{l+1,j}
 %&\\&=\sum^{n}_{j=p^l_{-}(A)}a_{l+1,j}-\sum^{n}_{j=p^l_{-}(A)+1}a_{l,j}-(\sum^{n}_{j=i}a_{l+1,j}-\sum^{n}_{j=i}a_{l,j})&\\&
=\varphi_l(A)-\langle \alpha_l^\vee,\wt(A)\rangle.\end{flalign*} Thus the map $\epsilon_l$ is already determined by $\varphi_l$ and conversely as well.
\endproof
\end{lem}
For the purpose of constructing an object in the category of crystals we define the Kashiwara operators by the following rule:
let $A$ be an arbitrary element of $B^{m,i}$ filled as in Definition~\ref{maindef}, then $\tilde{f}_{l}A$ and  $\tilde{e}_{l}A$ respectively is defined to be 0 if $\varphi_{l}(A)=0$ and $\epsilon_{l}(A)=0$ respectively. Otherwise the image of $A$ under $\tilde{f}_{l}$ and  $\tilde{e}_{l}$ respectively arises from $A$ by replacing certain boxes, namely
\begin{equation}\label{kashopf}
\tilde{f}_{l}A=\left\{
\begin{array}{ll}
\mbox{replace $\ytableausetup{boxsize=1.5em}\begin{ytableau}
a_{i,i}\end{ytableau}$ by \fbox{$a_{i,i}+1$}}, \quad \mbox{if $l=i$}\\
\mbox{replace \fbox{$a_{p^l_+(A),l-1}$} by 
\fbox{$a_{p^l_+(A),l-1}-1$} and 
\fbox{$a_{p^l_+(A),l}$} by 
\fbox{$a_{p^l_+(A),l}+1$}},\quad \mbox{if $l>i$}\\
\mbox{replace \fbox{$a_{l,p^l_-(A)}$} by \fbox{$a_{l,p^l_-(A)}+1$} and 
\fbox{$a_{l+1,p^l_-(A)}$} by \fbox{$a_{l+1,p^l_-(A)}-1$}}, \quad \mbox{if $l<i$}\\
\end{array}\right.
\end{equation}

\begin{equation}\label{kashope}
\tilde{e}_{l}A=\left\{
\begin{array}{ll}
\mbox{replace $\ytableausetup{boxsize=1.5em}\begin{ytableau}
a_{i,i}\end{ytableau}$ by \fbox{$a_{i,i}-1$}}, \quad \mbox{if $l=i$}\\
\mbox{replace \fbox{$a_{q^l_+(A),l-1}$} by \fbox{$a_{q^l_+(A),l-1}+1$} and 
\fbox{$a_{q^l_+(A),l}$} by \fbox{$a_{q^l_+(A),l}-1$}},\quad \mbox{if $l>i$}\\
\mbox{replace \fbox{$a_{l,q^l_-(A)}$} by \fbox{$a_{l,q^l_-(A)}-1$} and 
\fbox{$a_{l+1,q^l_-(A)}$}by \fbox{$a_{l+1,q^l_-(A)}+1$}}, \quad \mbox{if $l<i$}.\\

\end{array}\right.
\end{equation}

To be more accurate we should denote the Kashiwara operators by $_m\tilde{f}_l$ and $_m\tilde{e}_l$ respectively. However almost all Kashiwara operators, except $_m\tilde{f}_i$, are by the next lemma independent of $m$. Therefore the notation $\tilde{f}_l$ and $\tilde{e}_l$ respectively is justified. If there is no confusion we will also denote $_m\tilde{f}_i$ by $\tilde{f}_i$.
\begin{lem}\label{index}
Let $m,d\in\Z_{+}$ and $A\in B^{m,i}\cap B^{d,i},$ then we have  
$$_m\tilde{f}_lA=\ _d\tilde{f}_lA \mbox{ and } _m\varphi_l(A)=\ _d\varphi_l(A) \mbox{ for all $l\neq i$ },$$
$$_m\tilde{e}_lA=\ _d\tilde{e}_lA \mbox{ and } _m\epsilon_l(A)=\ _d\epsilon_l(A) \mbox{ for all $l$ }.$$
\proof
A short investigation of (\ref{ttt1}),(\ref{ttt2}),(\ref{kashopf}) and (\ref{kashope}) shows that they depend only on the filling of $A$ with the exception of $\varphi_i$ and $\tilde{f}_i$.
\endproof 
\end{lem}
\begin{rem}
\item It is not clear, why these operators are well-defined. Particulary we shall show in the next lemma that the images are always contained in $B^{m,i}$. 
\end{rem}
\begin{lem}\label{welldef}
For all $l\in I$ and $A\in B^{m,i}$ we have $\tilde{f}_{l}A,\tilde{e}_{l}A\in B^{m,i}$.
\proof
Assume that the element $\tilde{f}_{l}A$ and $\tilde{e}_{l}A$ respectively is not contained in $B^{m,i}$, then by definition there exists a Dyck path $(\beta(1),\dots, \beta(k))$, such that \begin{equation}\label{brag}\sum^k_{s=1}a_{\beta(s)}> m.\end{equation} This would be an impossible inequality if $l=i$; therefore we suppose that $l>i$, since the proof for $l<i$ is similar. By an inspection of the action of the Kashiwara operator $\tilde{f}_{l}$ we can conclude directly that (\ref{brag}) must be of the following form:
\begin{equation}\label{brag1}\sum^k_{s=1}a_{\beta(s)}=\sum^t_{s=1}a_{\beta(s)}+a_{z,l-1}+\sum^{p^l_+(A)}_{j=z}a_{j,l}+1+\sum^k_{s=t+p^l_+(A)-z+3}a_{\beta(s)},\end{equation} 
with an integer $z\in\{1,\cdots,i\}$ strictly smaller than $p^l_+(A)$ and some $t\in \{1,\cdots,k\}$. We get $$\sum^k_{s=1}a_{\beta(s)} > m\geq \sum^t_{s=1}a_{\beta(s)}+\sum^{p^l_+(A)}_{j=z}a_{j,l-1}+a_{p^l_+(A),l}+\sum^k_{s=t+p^l_+(A)-z+3}a_{\beta(s)}$$
and consequently
$$\sum^z_{j=1}a_{j,l-1}+\sum^{i}_{j=z}a_{j,l}=\sum^{p^l_+(A)}_{j=1}a_{j,l-1}+\sum^{i}_{j=p^l_+(A)}a_{j,l},$$ which is a contradiction to the choice of $p^l_+(A)$.
An inspection of the action with $\tilde{e}_{l}$ requires that (\ref{brag}) must be of the form:
$$\sum^k_{s=1}a_{\beta(s)}=\sum^t_{s=1}a_{\beta(s)}+\sum^{z}_{j=q^l_+(A)}a_{j,l-1}+1+a_{z,l}+\sum^k_{s=t+z-q^l_+(A)+3}a_{\beta(s)},$$
with an integer $z\in\{1,\cdots,i\}$ strictly greater than $q^l_+(A)$. Hence, together with
$$\sum^k_{s=1}a_{\beta(s)} > m\geq \sum^t_{s=1}a_{\beta(s)}+\sum^{z}_{j=q^l_+(A)}a_{j,l}+a_{q^l_+(A),l-1}+\sum^k_{s=t+z-q^l_+(A)+3}a_{\beta(s)},$$
we have again a contradiction to the choice of $q^l_+(A)$, namely
$$\sum^z_{j=1}a_{j,l-1}+\sum^{i}_{j=z}a_{j,l}=\sum^{q^l_+(A)}_{j=1}a_{j,l-1}+\sum^{i}_{j=q^l_+(A)}a_{j,l}.$$
\endproof
\end{lem}

Consequently, we have several well-defined maps which we need so as to prove our main result of Section~\ref{section3}. Before we state our theorem, we proof the following helpful lemma:

\begin{lem}\label{semir}
Let $A$ be an element in $B^{m,i}$, then we have $$\epsilon_l(A)=\max\{k\geq0|\tilde{e}_{l}^{k}A\neq 0\},\quad \varphi_l(A)=\max\{k\geq0|\tilde{f}_{l}^{k}A\neq 0\}.$$ 
\proof
As usual we proof only $\varphi_l(A)=\max\{k\geq0|\tilde{f}_{l}^{k}A\neq 0\}$ for $l>i$ because $l=i$ is trivial and the proof in the other cases are very similar. We will proceed by induction on $p^l_+(A)$. If $p^l_+(A)=1$ the proof is obvious so assume that $p^l_+(A)>1$ and let $r:=\min\{w\in\Z_{+} |p^l_+(\tilde{f}_{l}^{w}A)< p^l_+(A)\}.$ Then we obtain by the definition of $p^l_+(A)$ on the one hand $$r-1<\sum^{p^l_+(A)}_{j=p^l_+(\tilde{f}_{l}^{r}A)+1}a_{j,l-1}-\sum^{p^l_+(A)-1}_{j=p^l_+(\tilde{f}_{l}^{r}A)}a_{j,l}$$
and on the other hand, using the definition of $p^l_+(\tilde{f}_{l}^{r}A)$, we get $$r\geq\sum^{p^l_+(A)}_{j=p^l_+(\tilde{f}_{l}^{r}A)+1}a_{j,l-1}-\sum^{p^l_+(A)-1}_{j=p^l_+(\tilde{f}_{l}^{r}A)}a_{j,l}.$$
Hence the above inequality is actually a equality and by the induction hypothesis we can conclude 
$$\max\{k\geq0|\tilde{f}_{l}^{k}A\neq 0\}=r+\varphi_l(\tilde{f}^r_lA)=r+\sum^{p^l_+(\tilde{f}_{l}^{r}A)}_{j=1}a_{j,l-1}-\sum^{p^l_+(\tilde{f}_{l}^{r}A)-1}_{j=1}a_{j,l}=\varphi_l(A).$$
\endproof
\end{lem}

Now we are in position to state and to proof one of our main results in this paper, namely:
\begin{thm}\label{mainthm1}
The polytope $B^{m,i}$ together with the maps given by (\ref{weight}), (\ref{ttt1}), (\ref{ttt2}), (\ref{kashopf}) and (\ref{kashope}) becomes an abstract semiregular crystal.
\proof
The idea of the proof is to check step by step the properties (1)-(7) described in Definition~\ref{abcrystal}, whereby (2),(3) and (7) are obvious and (1),(4),(5) and the semiregularity are obvious with Lemma~\ref{fest} and Lemma~\ref{semir}. Thus it remains to prove the correctness of condition (6), whereby we verify as usual the statement only for $l>i$:
\vskip2pt
%$\bullet$ $\varphi_l(A)-\epsilon_l(A)=\langle h_l,wt(A)\rangle$\\
%This is an easy calculation since the right hand side is exactly
%$$\langle h_l,wt(A)\rangle=\begin{cases}
%m-\sum^i_{j=1}a_{j,i}-\sum^n_{j=i}a_{i,j}, & \text{if $l=i$}\\
%\sum^i_{j=1}a_{j,l-1}-\sum^i_{j=1}a_{j,l},& \text{if $l>i$}\\
%\sum^n_{j=i}a_{l+1,j}-\sum^n_{j=i}a_{l,i},& \text{if $l<i$}.
%\end{cases}$$
$\bullet$ $\tilde{f}_{l}A=A'$ if and only if $\tilde{e}_{l}A'=A$ for $A,A'\in B^{m,i}$\vskip2pt
Let $p^l_+(A)$ as in (\ref{1}) and let $q^l_+(\tilde{f}_{l}A)$ as in (\ref{2}). The assumption $q^l_+(\tilde{f}_{l}A)>p^l_+(A)$ gives $$\sum^{q^l_+(\tilde{f}_{l}A)}_{j=1}a_{j,l-1}+\sum^{i}_{j=q^l_+(\tilde{f}_{l}A)}a_{j,l}>\sum^{p^l_+(A)}_{j=1}a_{j,l-1}+\sum^{i}_{j=p^l_+(A)}a_{j,l},$$ which is a contradiction to the maximality.
The assumption $q^l_+(\tilde{f}_{l}A)<p^l_+(A)$ gives $$\sum^{q^l_+(\tilde{f}_{l}A)}_{j=1}a_{j,l-1}+\sum^{i}_{j=q^l_+(\tilde{f}_{l}A)}a_{j,l}\geq\sum^{p^l_+(A)}_{j=1}a_{j,l-1}+\sum^{i}_{j=p^l_+(A)}a_{j,l},$$ which is a contradiction to the minimality of $p^l_+(A).$ Now suppose similar as above that $p^l_+(\tilde{e}_{l}A)>q^l_+(A)$, then $$\sum^{p^l_+(\tilde{e}_{l}A)}_{j=1}a_{j,l-1}+\sum^{i}_{j=p^l_+(\tilde{e}_{l}A)}a_{j,l}\geq\sum^{q^l_+(A)}_{j=1}a_{j,l-1}+\sum^{i}_{j=q^l_+(A)}a_{j,l},$$ which is a contradiction to the maximality of $q^l_+(A)$ and $p^l_+(\tilde{e}_{l}A)<q^l_+(A)$ provides analogously $$\sum^{p^l_+(\tilde{e}_{l}A)}_{j=1}a_{j,l-1}+\sum^{i}_{j=p^l_+(\tilde{e}_{l}A)}a_{j,l}>\sum^{q^l_+(A)}_{j=1}a_{j,l-1}+\sum^{i}_{j=q^l_+(A)}a_{j,l},$$ which is a contradiction to the maximality.
Hence $p^l_+(A)=q^l_+(\tilde{f}_{l}A)$ and $q^l_+(A)=p^l_+(\tilde{e}_{l}A)$, which proves the theorem.
\endproof
\end{thm}

\begin{cor}\label{connected}
The crystal $B^{m,i}$ is connected.
\proof
It is immediate that for $A\in B^{m,i}$ with $\tilde{e}_{l}A=0$ for all $l\in I$, we must have $a_{p,q}=0$ for all $1\leq p \leq i\leq q \leq n$. Hence, for arbitrary elements $A$ and $B$ there exists always a couloured path from $A$ to $B$.
\endproof
\end{cor}
%%%%%%%%%%%%%%%%%%%%%%%%%%%%%%%%%%%%%%%%%%%%%%%%%%%%%%%%%%%%%%%%%%%%%%%%%%%%%%%%%%%%%%%%%%%%%%%%%%%%%%%%%%%%%%%%%%%%%%%%%%%%%%%%%%%
%         
%%%%%%%%%%%%%%%%%%%%%%%%%%%%%%%%%%%%%%%%%%%%%%%%%%%%%%%%%%%%%%%%%%%%%%%%%%%%%%%%%%%%%%%%%%%%%%%%%%%%%%%%%%%%%%%%%%%%%%%%%%%%%%%%%%%
\section{Tensor products and Nakajima monomials}\label{section4}
In this section, we want to recall tensor products of crystals and investigate the action of Kashiwara operators on tensor products. Furthermore, we want to introduce a crystal, the set of all Nakajima monomials, such that we can think of $B(\lambda)$, where $\lambda$ is a dominant integral $A_n$ weight, as a set of certain monomials. This theory is discovered by Nakajima \cite{Nak03}, and generalized by Kashiwara \cite{Kashi03} and will be important in the following sections.
\subsection{Tensor product of crystals}
Suppose that we have two abstract crystals $B_1$, $B_2$ in the sense of Definition~\ref{abcrystal}, then we can construct a new crystal which is as a set nothing but $B_1\times B_2$. This crystal is denoted by $B_1\otimes B_2$ and the Kashiwara operators are given as follows:
$$\tilde{f}_{l}(b_1\otimes b_2)=\begin{cases} (\tilde{f}_{l}b_1)\otimes b_2, \text{ if $\varphi_l(b_1)>\epsilon_l(b_2)$}\\
b_1\otimes (\tilde{f}_{l}b_2), \text{ if $\varphi_l(b_1)\leq \epsilon_l(b_2)$}\end{cases}$$
$$\tilde{e}_{l}(b_1\otimes b_2)=\begin{cases} (\tilde{e}_{l}b_1)\otimes b_2, \text{ if $\varphi_l(b_1)\geq \epsilon_l(b_2)$}\\
b_1\otimes (\tilde{e}_{l}b_2), \text{ if $\varphi_l(b_1)<\epsilon_l(b_2)$.}\end{cases}$$
Further, one can describe explicitely the maps $\wt$,$\varphi_l$ and $\epsilon_l$ on $B_1\otimes B_2$, namely:
$$\wt(b_1\otimes b_2)=\wt(b_1)+\wt(b_2)$$
$$\varphi_l(b_1\otimes b_2)=\max\{\varphi_l(b_2),\varphi_l(b_1)+\varphi_l(b_2)-\epsilon_l(b_2)\}$$
$$\epsilon_l(b_1\otimes b_2)=\max\{\epsilon_l(b_1),\epsilon_l(b_1)+\epsilon_l(b_2)-\varphi_l(b_1)\}.$$
A very important point in representation theory is to determine crystal bases for irreducible modules over quantum algebras. This leads to many  combinatorial models, discovered in a series of papers (\cite{NK94},\cite{L94}, \cite{NS97}). Since this paper has the goal of determining the crystal graph of KR modules, we will only mention in the following remark how we could compute crystal bases for $V(\lambda)$ using tensor products of polytopes.
\begin{rem}
Remember that we have already a crystal structure on the sets $B^{m,i}$ for all $i=1,\cdots,n$ and by the above considerations also on $$B^{m_1,1}\otimes\cdots\otimes B^{m_n,n}.$$ It means we are considering patterns
\vskip 8pt
$\yng(1,3,6,10,15,21,28)$
\vskip 10pt

with the following crystal structure:
let us take such a pattern $\mathbf A:=(A_1,\cdots,A_n)$ and fix $l\in I$, then we assign to $\mathbf A$ a sequence of - `s followed by a sequence of +`s
$$\mbox{seq($\mathbf A$)}:=(\underbrace{-,\cdots,-}_{\epsilon_l(A_1)},\underbrace{+,\cdots,+}_{\varphi_l(A_1)},\underbrace{-,\cdots,-}_{\epsilon_l(A_2)},\underbrace{+,\cdots,+}_{\varphi_l(A_2)},\cdots,\underbrace{-,\cdots,-}_{\epsilon_l(A_n)},\underbrace{+,\cdots,+}_{\varphi_l(A_n)})$$
and cancel out all $(+,-)$-pairs to obtain the so called $\ell$-signature
\begin{equation}\label{signature}\mbox{$\ell$-sgn($\mathbf A$)}:=(-,\cdots,-,+,\cdots,+).\end{equation}
Using the $\ell$-signature \mbox{$\ell$-sgn($\mathbf A$)} of $\mathbf A$, we can describe the Kashiwara operators. Assume the left-most $+$ in the $\ell$-signature corresponds to $A_i$ and the right-most $-$ corresponds to $A_j$, then 
$$\tilde{f}_{l}(A_1A_2\ldots A_n)=A_1\ldots A_{i-1}(\tilde{f}_{l}A_i)A_{i+1}\ldots A_n$$
$$\tilde{e}_{l}(A_1A_2\ldots A_n)=A_1\ldots A_{j-1}(\tilde{e}_{l}A_j)A_{j+1}\ldots A_n.$$
Hence we have defined a new crystal, and with Theorem~\ref{nakmon} it is clear from standard arguments that the connected component of $0$ (all boxes filled with 0) is isomorphic to the crystal $B(\lambda)$.
\end{rem}

\subsection{Nakajima monomials}

For $i\in I$ and $n\in\Z$ we consider monomials in the variables $Y_i(n)$, i.e. we obtain the set of Nakajima monomials $\mathcal{M}$ as follows:
$$\mathcal{M}:=\{\prod_{i\in I, n\in\Z} Y_i(n)^{y_i(n)}| y_i(n)\in\Z \mbox{ vanish except for finitely many $(i,n)$}\}$$
With the goal to define the crystal structure on $\mathcal{M}$, we take some integers $c=(c_{i,j})_{i\neq j}$ such that $c_{i,j}+c_{j,i}=1$. Let now $M=\prod_{i\in I, n\in\Z} Y_i(n)^{y_i(n)}$ be an arbitrary monomial in $\mathcal{M}$ and $l\in I$, then we set:
$$\wt(M)=\sum_i(\sum_n y_i(n))\omega_i$$
$$\varphi_l(M)=\max\{\sum_{k\leq n}y_l(k)|n\in\Z\}$$
$$\epsilon_l(M)=\max\{-\sum_{k>n}y_l(k)|n\in\Z\}$$
and 
$$n^l_f=\min\{n|\varphi_l(M)=\sum_{k\leq n}y_l(k)\}$$
$$n^l_e=\max\{n|\epsilon_l(M)=-\sum_{k>n}y_l(k)\}$$
The Kashiwara operators are defined as follows:
$$\tilde{f}_{l}M=\begin{cases}
A_l(n^l_f)^{-1}M, \text{ if $\varphi_l(M)>0$}\\
0, \text{ if $\varphi_l(M)=0$}\end{cases}$$

$$\tilde{e}_{l}M=\begin{cases}
A_l(n^l_e)M, \text{ if $\epsilon_l(M)>0$}\\
0, \text{ if $\epsilon_l(M)=0$,}\end{cases}$$
whereby $$A_l(n):=Y_l(n)Y_l(n+1)\prod_{i\neq l}Y_i(n+c_{i,l})^{\langle \alpha_i^{\vee},\alpha_l\rangle}.$$
The following two results due to Kashiwara \cite{Kashi03} are essential for the process of this paper:
\begin{prop}
With the maps $\wt$, $\varphi_l, \epsilon_l, \tilde{f}_{l},\tilde{e}_{l}$, $l\in I$, the set $\mathcal{M}$ becomes a semiregular crystal.
\end{prop}

\begin{rem}
\textit{A priori} the crystal structure depends on $c$, hence we will denote this crystal by $\mathcal{M}_c$. But it is easy to see that the isomorphism class of $\mathcal{M}_c$ does not depend on this choice. In the literature $c$ is often chosen as
$$c_{i,j}=\begin{cases} 0, \text{ if $i>j$}\\
 1, \text{ else} \end{cases} \quad or \quad c_{i,j}=\begin{cases} 0, \text{ if $i<j$}\\
 1, \text{ else.} \end{cases}$$
\end{rem}
\begin{prop}\label{iii}
Let $M$ be a monomial in $\mathcal{M}$, such that $\tilde{e}_{l}M=0$ for all $l\in I$. Then the connected component of $\mathcal{M}$ containing $M$ is isomorphic to $B(\wt(M))$.
\end{prop}
%%%%%%%%%%%%%%%%%%%%%%%%%%%%%%%%%%%%%%%%%%%%%%%%%%%%%%%%%%%%%%%%%%%%%%%%%%%%%%%%%%%%%%%%%%%%%%%%%%%%%%%%%%%%%%%%%%%%%%%%%%%%%%%%%%%
%         
%%%%%%%%%%%%%%%%%%%%%%%%%%%%%%%%%%%%%%%%%%%%%%%%%%%%%%%%%%%%%%%%%%%%%%%%%%%%%%%%%%%%%%%%%%%%%%%%%%%%%%%%%%%%%%%%%%%%%%%%%%%%%%%%%%%
\section{Stembridge axioms and isomorphism of crystals}\label{section5}
By using the description of crystal graphs by certain monomials and maps on these, we want to show that our set $B^{m,i}$ satisfies the so called Stembridge axioms stated in \cite{S2003}. These axioms give a local characterization of simply-laced crystals which are helpful if one could not find a global isomorphism. Since we could not find an isomorphism from $B^{m,i}$ to the set of all monomials describing $B(m\omega_i)$, we will identify only certain $A_2$-crystals. Let us first recall a slightly modified result from \cite{S2003}.
\subsection{Stembridge axioms}
The basic idea of the following proposition is to give a simple set of local axioms to characterize the set of crystals of representations in the class of all crystals. In particular, with these axioms one can determine whether or not a crystal is the crystal of a representation. 
\begin{prop}
Let $\Lg$ be a simply-laced Lie algebra and $B$ be a connected crystal graph, such that the following conditions are satisfied:
\begin{enumerate}
\item If $\tilde{e}_{l}b$ is defined, then $\epsilon_j(\tilde{e}_{l}b)\geq \epsilon_j(b)$ and $\varphi_j(\tilde{e}_{l}b)\leq \varphi_j(b)$ for all $j\neq l$.
\item If $\tilde{e}_{l},\tilde{e}_{j}b$ are defined and $\epsilon_j(\tilde{e}_{l}b)=\epsilon_j(b)$, then $\tilde{e}_{l}\tilde{e}_{j}b=\tilde{e}_{j}\tilde{e}_{l}b$ and $\varphi_l(b^{'})=\varphi_l(\tilde{f}_{j}b^{'}),$ where $b^{'}=\tilde{e}_{l}\tilde{e}_{j}b=\tilde{e}_{j}\tilde{e}_{l}b$.
\item If $\tilde{e}_{l},\tilde{e}_{j}b$ are defined and $\epsilon_j(b)-\epsilon_j(\tilde{e}_{l}b)=\epsilon_l(b)-\epsilon_l(\tilde{e}_{j}b)=-1$, then $\tilde{e}_{l}\tilde{e}_{j}^2\tilde{e}_{l}b=\tilde{e}_{j}\tilde{e}_{l}^2\tilde{e}_{j}b$ and $\varphi_l(b^{'})-\varphi_l(\tilde{f}_{j}b^{'})=\varphi_j(b^{'})-\varphi_j(\tilde{f}_{l}b^{'})=-1,$ where $b^{'}=\tilde{e}_{l}\tilde{e}_{j}^2\tilde{e}_{l}b=\tilde{e}_{j}\tilde{e}_{l}^2\tilde{e}_{j}b$.
\item If $\tilde{f}_{l},\tilde{f}_{j}b$ are defined and $\varphi_j(\tilde{f}_{l}b)=\varphi_j(b)$, then $\tilde{f}_{l}\tilde{f}_{j}b=\tilde{f}_{j}\tilde{f}_{l}b$ and $\epsilon_l(b^{'})=\epsilon_l(\tilde{e}_{j}b^{'}),$ where $b^{'}=\tilde{f}_{l}\tilde{f}_{j}b=\tilde{f}_{j}\tilde{f}_{l}b$.
\item If $\tilde{f}_{l},\tilde{f}_{j}b$ are defined and $\varphi_j(b)-\varphi_j(\tilde{f}_{l}b)=\varphi_l(b)-\varphi_l(\tilde{f}_{j}b)=-1$, then $\tilde{f}_{l}\tilde{f}_{j}^2\tilde{f}_{l}b=\tilde{f}_{j}\tilde{f}_{l}^2\tilde{f}_{j}b$ and $\epsilon_l(b^{'})-\epsilon_l(\tilde{e}_{j}b^{'})=\epsilon_j(b^{'})-\epsilon_j(\tilde{e}_{l}b^{'})=-1,$ where $b^{'}=\tilde{f}_{l}\tilde{f}_{j}^2\tilde{f}_{l}b=\tilde{f}_{j}\tilde{f}_{l}^2\tilde{f}_{j}b$.
\end{enumerate}
Then B is a crystal graph induced by a representation.
\end{prop}

\subsection{Isomorphism of \texorpdfstring{$A_2$}{A2} crystals}
We want to make full use of the above mentioned result to prove the following main theorem of this section:

\begin{thm}\label{nakmon}
We have an isomorphism of crystals $$B^{m,i}\cong B(m\omega_i).$$
\end{thm}
\proof
Assume $|j-r|=1$, i.e. $r=j+1$ and let $A$ be an arbitrary element in $B^{m,i}$. Since $B^{m,i}$ is a connected crystal we cancel all arrows with colour $s\neq j,j+1$ and denote the remaining $(j,j+1)$-connected graph containing $A$ by $Z_{(j,j+1)}(A)$. We define a map $\Psi: Z_{(j,j+1)}(A)\cup \{0\}\longrightarrow \mathcal{M}\cup \{0\}$ by mapping $0$ to $0$ and $B$ to:

$$ \begin{cases} Y_1(i)^{m-\sum^n_{s=i}b_{i,n}}\prod^i_{k=1}Y_2(k)^{b_{k,i}}Y_2(k+1)^{-b_{k,i+1}}\prod^i_{k=1}Y_1(k+1)^{-b_{k,i}}, \text{ if $j=i$}\\\\
Y_2(n-1)^{m-\sum^i_{s=1}b_{s,i}}\prod^n_{k=i}Y_1(k)^{b_{i,n+i-k}}Y_1(k+1)^{-b_{i-1,n+i-k}}\prod^n_{k=i}Y_2(k)^{-b_{i,n+i-k}}, \text{ if $j=i-1$}\\\\
\prod^i_{k=1}Y_1(k)^{b_{k,j-1}}Y_1(k+1)^{-b_{k,j}}\prod^i_{k=1}Y_2(k)^{b_{k,j}}Y_2(k+1)^{-b_{k,j+1}}, \text{ if $j>i$}\\\\
\prod^n_{k=i}Y_1(k)^{b_{j+1,n+i-k}}Y_1(k+1)^{-b_{j,n+i-k}}\prod^n_{k=i}Y_2(k)^{b_{j+2,n+i-k}}Y_2(k+1)^{-b_{j+1,n+i-k}}, \text{ if $j+1<i$.}\end{cases}$$

We are claiming that $\Psi$ is an $A_2\cong \mathfrak{sl}_{3}(j,j+1)$ crystal isomorphism. Note that $\Psi$ has the following properties:\vskip 5pt

$\bullet$ $\wt(\Psi(B))=\wt(B)$ \vskip 3pt
$\bullet$ $\varphi_l(\Psi(B))=\varphi_l(B), \epsilon_l(\Psi(B))=\epsilon_l(B)$. The proof is a case-by-case consideration, for instance if $l=j>i$, then $$\varphi_l(\Psi(B))=\max\{\sum^r_{s=1}b_{s,j-1}-\sum^{r-1}_{s=1}b_{s,j}|1\leq r \leq i\}=\varphi_l(B),$$ because 
the maximum occurs at least at $r=p_{+}^l(B)$. \vskip 3pt
$\bullet$ $\Psi$ commutes with the Kashiwara operators: by Lemma~\ref{semir} and the above computations we can conclude that $\Psi$ commutes with all $\tilde{f}_{l},\tilde{e}_{l}$ acting by zero on $B$. So assume $\tilde{f}_{l}B=\hat{B}$ and $\tilde{e}_{l}B=\hat{B}$ respectively. Our aim is to prove $\Psi(\hat{B})=A_l(n^l_f)^{-1}\Psi(B)$, $\Psi(\hat{B})=A_l(n^l_e)\Psi(B)$. This is again a case-by-case consideration, for instance if $l=j+1<i$ we have $n^l_f=n+i-p_{-}^l(B)$ and

\begin{flalign*} &\Psi(\hat{B})=\prod^n_{k=i}Y_1(k)^{\hat{b}_{j+1,n+i-k}}Y_1(k+1)^{-\hat{b}_{j,n+i-k}}\prod^n_{k=i}Y_2(k)^{\hat{b}_{j+2,n+i-k}}Y_2(k+1)^{-\hat{b}_{j+1,n+i-k}}
&\\&=\prod^n_{k=i, k\neq n^l_f}Y_1(k)^{b_{j+1,n+i-k}}\prod^n_{k=i}Y_1(k+1)^{-b_{j,n+i-k}} \prod^n_{k=i, k\neq n^l_f}Y_2(k)^{b_{j+2,n+i-k}}Y_2(k+1)^{-b_{j+1,n+i-k}}&\\& \times Y_1(n^l_f)^{b_{j+1,p_-^l(B)}+1}Y_2(n^l_f)^{b_{j+2,p_-^l(B)}-1}Y_2(n^l_f+1)^{-b_{j+1,p_-^l(B)}-1}&\\&=Y_1(n^l_f)Y_2(n^l_f)^{-1}Y_2(n^l_f+1)^{-1}\Psi(B)=A_l(n^l_f)^{-1}\Psi(B).\end{flalign*}
Hence $\Psi$ is a strict crystal morphism.\vskip 3pt
$\bullet$ $\Psi$ is bijective: since $Z_{(j,j+1)}(A)$ is connected and $\Psi$ is a crystal morphism we get that $\mathcal{I}m(\Psi)$ is connected and contains at least, and therefore by Proposition~\ref{iii} one highest weight monomial, say of weight $\mu$. So the image is isomorphic, again by Proposition~\ref{iii} to the $\mathfrak{sl}_{3}(j,j+1)$ crystal $B(\mu)$. Let $T\in Z_{(j,j+1)}(A)$ be a highest weight element, such that $\tilde{e}_{i_1}\cdots \tilde{e}_{i_s}A=T$, then the restriction of $\Psi$ to $G=\{\tilde{f}_{j_1}\cdots \tilde{f}_{j_s}T|j_1,\cdots,j_s\in\{j,j+1\}\}$ is an isomorphism. However, we have $Z_{(j,j+1)}(A)=G$.\vskip 5pt
Therefore, we can conclude that $A$ satisfies the Stembridge axioms for all $j,r\in I$ with $|j-r|=1$, whereby the other relations are easily verified.
\endproof
Summerizing we have defined a set $B^{m,i}$ which is by Theorem~\ref{mainthm1} a crystal and by Theorem~\ref{nakmon} actually the crystal $B(m\omega_i)$. In the following we want to collect some known facts about KR-crystals and define the Kashiwara operators $\tilde{f}_0,\tilde{e}_0$ on our underlying polytope.
%%%%%%%%%%%%%%%%%%%%%%%%%%%%%%%%%%%%%%%%%%%%%%%%%%%%%%%%%%%%%%%%%%%%%%%%%%%%%%%%%%%%%%%%%%%%%%%%%%%%%%%%%%%%%%%%%%%%%%%%%%%%%%%%%%%
%         
%%%%%%%%%%%%%%%%%%%%%%%%%%%%%%%%%%%%%%%%%%%%%%%%%%%%%%%%%%%%%%%%%%%%%%%%%%%%%%%%%%%%%%%%%%%%%%%%%%%%%%%%%%%%%%%%%%%%%%%%%%%%%%%%%%%
\section{The promotion operator}\label{section6}
The existence of KR-crystals of type $A^{(1)}_n$ was shown in \cite{KKMMNN92} and a combinatorial description was provided in \cite{S02},\cite{KW12}. The existence of KR-crystals for non-exceptional types can be found in \cite{O07},\cite{OS08} and further a combinatorial description is provided in \cite{FOS09}. Summerizing the results for type $A^{(1)}_n$, a model for KR-crystals is given by the set of all semi-standard Young tableuax of shape $\lambda=m\omega_i$ with affine Kashiwara operators \begin{equation}\label{affka}\tilde{f}_{0}:=pr^{-1}\circ \tilde{f}_{1}\circ pr, \mbox{ and } \tilde{e}_{0}:=pr^{-1}\circ \tilde{e}_{1}\circ pr,\end{equation} whereby $pr$ is the so called Schützenbergers's promotion operator \cite{S72}, which is the analogue of the cyclic Dynkin diagram automorphism on the level of crystals. The promotion operator on the set of all semi-standard Young tableaux over the alphabet $1\prec 2\cdots \prec n+1$ can be obtained by using jeu-de-taquin. Particulary let $T$ be a Young tableaux, then we get $pr(T)$ by removing all letters $n+1$, adding 1 to each letter in the remaining tableaux, using jeu-de-taquin to slide all letters up and finally filling the holes with 1's. Our aim now is to define the Schützenberger promotion operator on our polytope $B^{m,i}$, to obtain a polytope realization of these crystals. Before we are in position to define such a map we will first state an important result due to \cite{BST10}. For simplicity we write $\wt(A)= (r_1,\cdots,r_{n+1})$, if $\wt(A)=\sum^{n+1}_{j=1}r_j\epsilon_j$, whereby $\epsilon_j:\mathfrak{sl}_{n+1}\longrightarrow \C$ is the projection on the $j$-th diagonal entry.
\begin{prop}\label{prch}
Let $\Psi:B(m\omega_i)\longrightarrow B(m\omega_i)$ be a map, such that
\begin{enumerate}
\item $\Psi$ shifts the content, which means if $\wt(A)=(r_1,\cdots,r_{n+1})$, then $\wt(\psi(A))=(r_{n+1},r_1,\cdots,r_{n})$
\item $\Psi$ is bijective
\item $\Psi\circ\tilde{f}_{j}=\tilde{f}_{j+1}\circ \Psi$, $\Psi\circ\tilde{e}_{j}=\tilde{e}_{j+1}\circ \Psi$ for all $j\in\{1,\cdots,n-1\},$
\end{enumerate}
then $\Psi=pr$.
\end{prop}
In particular, using Theorem~\ref{nakmon}, this means that it is sufficient to define a map on $B^{m,i}$ satisfying the conditions (1)-(3). The computation of this map, which we will denote already by $pr$, will proceed by an algorithm consisting of $i$ steps, where each step will give us a column of $pr(A)$. We denote by $pr(a_{r,s})$ the entries of $pr(A)$ and by $\mathbf a_j$ and $pr(\mathbf{a_j})$ respectively the $j$-th column of $A$ and $j$-th column of $pr(A)$ respectively for $j=1,\cdots,i$ and mostly we assume the notation explained in Remark~\ref{notrem}. In order to state the algorithm we denote further by $(\mathbf a_j)^{\geq l}$ the column obtained from $\mathbf a_j$ by canceling all entries between $i$ and $l-1$. For instance if $\mathbf a_j$ is the $j$-th column of some $A\in B^{m,4}$ we have $$(\mathbf a_j)^{\geq 5}= \young(5,0),\mbox{  if  } a_j=\young(3,5,0).$$

If $i=1$, then it is obvious that the map $pr$ defined by $$pr(a_{1,j})=\begin{cases} m-\sum^n_{r=1}a_{1,r},& \text{ if $j=1$}\\
a_{1,j-1},& \text{ else}
\end{cases}$$
satisfies the conditions (1)-(3). So for $n\geq i\geq 2$ we will use the following algorithm to compute $pr$:
\vskip13pt

\subsection{Algorithm}\label{algorithm}
For given $A\in B^{m,i}$ we implement the following steps to compute $pr(A)$:
\begin{enumerate}[\bfseries(1)]
\item Consider the $(i-1)$-th and $i$-th column of $A$ and compute inductively the integers $i\leq l^{i-1}_1<l^{i-1}_2<\cdots<l^{i-1}_{t_{i-1}}=n$, whereby $$l^{i-1}_j=q_-((\mathbf a_{i-1})^{> l^{i-1}_{j-1}},(\mathbf a_{i})^{> l^{i-1}_{j-1}})$$ and where we undestand $l^{i-1}_0=i-1$. The $i$-th column of $pr(A)$ is then given by: $$pr(a_{i,r})=\begin{cases} \epsilon(\mathbf a_{i-1},\mathbf a_{i}),& \text{if $r=i$}\\ 
\epsilon((\mathbf a_{i-1})^{\geq r},(\mathbf a_{i})^{\geq r}),& \text{if $r-1\in\{l^{i-1}_1,\cdots,l^{i-1}_{t_{i-1}-1}\}$}\\ 
a_{i,r-1}, & \text{else. }
\end{cases}$$

Further define a new column 
\begin{equation}\label{newcolumn}\widehat{a_{i-1,r}}=\begin{cases} a_{i-1,r}+a_{i,r}-\epsilon((\mathbf a_{i-1})^{> r},(\mathbf a_{i})^{>r}),&  \text{ if $r\in\{l^{i-1}_1,\cdots,l^{i-1}_{t_{{i-1}}-1}\}$}\\ 
a_{i-1,n}+a_{i,n},& \text{if $r=n$}\\ 
a_{i-1,r}, & \text{else. }
\end{cases}\end{equation}

\item With the aim to determine the $(i-1)$-th column of $pr(A)$ repeat step (1) with the $(i-2)$-th column of $A$ and the new defined column (\ref{newcolumn}). Using the integers $i\leq l^{i-2}_1<l^{i-2}_2<\cdots<l^{i-2}_{t_{i-2}}=n$ compute the $(i-2)$-th column as in step (1).
\item Repeat step (2) as long as all columns, except the first one, from $pr(A)$ are known.
\item The first column of $pr(A)$ is given as follows:

$$pr(a_{1,r})=\begin{cases} m-\sum^n_{j=i}a_{i,j}-\sum^i_{j=2}pr(a_{j,i}),& \text{if $r=1$} \\ 
 \sum^i_{j=1}a_{j,r-1}-\sum^i_{j=2}pr(a_{j,r}),& \text{if $r>1$.}
\end{cases}$$
\end{enumerate}

\begin{rem}\label{imrem}\mbox{}
\begin{enumerate}
\item 
\textit{A priori} it is not clear, why the entries in the first column are non negative integers. This will be our first step and is proven in Proposition~\ref{wise}.
\item As well it is not clear, why the image of $pr$ lies in $B^{m,i}$. With the purpose to prove the well-definedness of $pr$ we will show first for any $A\in B^{m,i}$ that $pr \circ\tilde{e}_{j}A=\tilde{e}_{j+1}\circ pr(A)$ holds for $j=1,\cdots,n-1$, where the equation can be understood by Lemma~\ref{index} as a equation independently from the knowledge where $pr(A)$ lives.
\end{enumerate}
\end{rem}
\begin{ex}
\begin{enumerate}[i)]
\item
We pick one element $A$ (see below) in $B^{3,3}$ and follow our algorithm. In the first step we get $l^2_{1}=3<l^2_{2}=4<l^2_{3}=5$ which gives us the third column:
$$
{\Yvcentermath1
A=\young(111,200,000)\rightsquigarrow\young(\bullet\bullet 1,\bullet\bullet 0,\bullet\bullet 0)
}
$$
The new column is given by \young(2,0,0). So following step two we get again $l^1_{1}=3<l^1_{2}=4<l^1_{3}=5$ and hence
$$
{\Yvcentermath1
\young(111,200,000)\rightsquigarrow\young(\bullet\bullet 1,\bullet\bullet 0,\bullet\bullet 0)\rightsquigarrow\young(\bullet11,\bullet20,\bullet00)
}
$$
So our last step gives us $$pr(A)=\young(011,120,200)$$
Another example in $B^{4,3}$ is described below:
$${
\Yvcentermath1
\young(110,011,000)\rightsquigarrow\young(\bullet\bullet2,\bullet\bullet0,\bullet\bullet0)\rightsquigarrow\young(\bullet12,\bullet00,\bullet00)\rightsquigarrow\young(012,200,200)
}$$
\item 
Now we pick an element in $B^{7,4}$, namely $$\young(1011,0132,1020)$$ and get $l^3_{1}=5<l^3_{2}=6$, $l^2_{1}=4<l^2_{2}=5<l^2_{3}=6$ and $l^1_{1}=4<l^1_{2}=5<l^1_{3}=6$. The new columns after the first step and second step respectively are given by $$\young(1,3,2) \mbox{ and }  \young(0,4,2) \mbox{ respectively.}$$
Thus, step by step we obtain the columns of $pr(A)$:

$${
\Yvcentermath1
\young(1011,0132,1020)\rightsquigarrow\young(\bullet\bullet\bullet3,\bullet\bullet\bullet1,\bullet\bullet\bullet2)\rightsquigarrow\young(\bullet\bullet03,\bullet\bullet11,\bullet\bullet02)\rightsquigarrow\young(\bullet103,\bullet011,\bullet102)\rightsquigarrow\young(0103,1011,3102)
}$$
\end{enumerate}
\end{ex}

\begin{prop}\label{wise}
For any $A\in B^{m,i}$ the first column of $pr(A)$ consists of non-negative integers and moreover $$m-\sum^n_{r=i}pr(a_{1,r})=\sum^i_{r=1}a_{r,n}.$$
\proof
We will prove this statement by induction on $i$. Assume $i=2$ and let $l^{1}_1<l^{1}_2<\cdots<l^{1}_{t_{1}}$ be the integers decribed in the algortihm. The first entry in the first column is exactly $$m-\sum^n_{r=2}a_{2,r}-\epsilon(\mathbf a_{1},\mathbf a_{2})=m-\sum^{l^1_1}_{r=2}a_{1,r}-\sum^n_{r=l^1_1}a_{2,r}\geq 0,$$ and the other entries in the first column of $pr(A)$ are either of the form $a_{1,k}$ or of the form $a_{1,l^{1}_k}+a_{2,l^{1}_k}-\epsilon((\mathbf a_{1})^{>l^1_k},(\mathbf a_{2})^{>l^1_k})$, for some $k$. However, by the definition of $l^1_{k}$ we know that $a_{2,l^{1}_k}+\cdots+a_{2,l^{1}_{k+1}-1}\geq a_{1,l^{1}_k+1}+\cdots+a_{1,l^{1}_{k+1}}$, which gives us $a_{2,l^{1}_k}\geq\epsilon((\mathbf a_{1})^{>l^1_k},(\mathbf a_{2})^{>l^1_k})$. Thus all entries in the first column are non negative integers. Now we will show that for all $j\in\{0,\cdots,t_1-1\}$ \begin{flalign*}m-\sum^n_{r=2}pr(a_{1,r})=\sum^{l^1_{j+1}}_{r=l^1_{j}+1}a_{1,r}+\sum^n_{r=l^1_{j+1}}a_{2,r}-\sum^n_{r=l^1_{j}+2}pr(a_{1,r})\end{flalign*} holds and proceed by upward induction on $j$. The $j=0$ case is obvious and assuming that $j>0$ we obtain by using the induction hypothesis
\begin{flalign*}&m-\sum^n_{r=2}pr(a_{1,r})=\sum^{l^1_{j}}_{r=l^1_{j-1}+1}a_{1,r}+\sum^n_{r=l^1_{j}}a_{2,r}-\sum^n_{r=l^1_{j-1}+2}pr(a_{1,r})&\\&=\sum^{l^1_{j}}_{r=l^1_{j-1}+1}a_{1,r}+\sum^n_{r=l^1_{j}}a_{2,r}-\sum^{l^1_{j}+1}_{r=l^1_{j-1}+2}pr(a_{1,r})-\sum^n_{r=l^1_{j}+2}pr(a_{1,r})&\\&=a_{1,l^1_j}+\sum^n_{r=l^1_{j}}a_{2,r}-pr(a_{1,l^1_j+1})-\sum^n_{r=l^1_{j}+2}pr(a_{1,r})&\\&=\sum^{l^1_{j+1}}_{r=l^1_{j}+1}a_{1,r}+\sum^n_{r=l^1_{j+1}}a_{2,r}-\sum^n_{r=l^1_{j}+2}pr(a_{1,r}),\end{flalign*} which finishes the induction. According to this we complete the initial step, since $$m-\sum^n_{r=2}pr(a_{1,r})=\sum^{l^1_{t_1}}_{r=l^1_{t_1-1}+1}a_{1,r}+\sum^n_{r=l^1_{t_1}}a_{2,r}-\sum^n_{r=l^1_{t_1-1}+2}pr(a_{1,r})=a_{1,n}+a_{2,n}.$$
%We will prove in additional for completeness the well-definedness in the $i=2$ case. Let $(\beta(1),\dots, \beta(k))$ be a Dyck path, such that $\beta(1)=\alpha_{1,2},\cdots,\beta(z)=\alpha_{1,z+1}$,\ $l^1_{p-1}<z+1\leq l^1_p$, and the rest is of the form $\alpha_{2,\bullet}$ and where we understand $l^1_{0}=i-1=1$. If $z>l^1_{p-1}$ we have $\sum pr(a_{\beta(k)})\leq \sum pr(a_{\beta'(k)})$, whereby  $$\beta'(k)=\begin{cases} \beta(k),& \text{if $k\leq z$ or $k\geq l_p+1$} \\ 
%pr(\alpha_{1,k+1}),& \text{if $z<k\leq l_p$}
%\end{cases}$$
%Thereafter we assume that $z=l^1_{p-1}$ and get 
%\begin{flalign*}\sum pr(a_{\beta(k)})&=m-\sum^{n}_{w=l^1_{p}}a_{2,w}-\sum^{t_1-1}_{r=p}(\sum^{l^1_{r+1}}_{w=l^1_r+1}a_{1,w}-\sum^{l^1_{r+1}-1}_{w=l^1_r+1}a_{2,w})&\\& =m-\sum^{n}_{w=l^1_{p}}a_{2,w}+\sum^{n}_{w=l^1_{p}+1}a_{1,w}\leq m\end{flalign*} Hence we have $pr(A)\in B^{m,2}$. 
Now let $i>2$ and consider an element $B\in B^{m,i-1}$, constructed as follows. The first $(i-2)$-columns of $B$ are the same as the ones from $A$ and the $(i-1)$-th column is precisely the new obtained column in (\ref{newcolumn}), which we get if we apply step (1) of the algorithm to $A$. In other words, we erase the $i$-th column of $A$ and replace the $(i-1)$-th column of $A$ by $$b_{i-1,r}=\begin{cases} a_{i-1,r}+a_{i,r}-\epsilon((\mathbf a_{i-1})^{> r},(\mathbf a_{i})^{>r}),&  \text{ if $r\in\{l^{i-1}_1,\cdots,l^{i-1}_{t_{{i-1}}-1}\}$}\\ 
a_{i-1,n}+a_{i,n},& \text{if $r=n$}\\ 
a_{i-1,r}, & \text{else }
\end{cases}$$ and obtain $B$. One can easily check that $B\in B^{m,i-1}$, where $B^{m,i-1}$ is the polytope associated to the Lie algebra $A_{n-1}$. Particulary we claim the following: if we glue the $i$-th column of $pr(A)$ to $pr(B)$ the resulting element is again $pr(A)$, i.e. $$pr(A)=pr(B)pr(\mathbf{a_i}).$$By the definition of the algorithm the claim is obvious for all columns except the first one. Because of that let $\ytableausetup{boxsize=1.4em}
\begin{ytableau}
b_{i} & \dots
&  b_{n} \\
\end{ytableau}$ be the transpose of the first column of $pr(B)$ and $\ytableausetup{boxsize=1.4em}
\begin{ytableau}
a_{i} & \dots
& a_{n} \\
\end{ytableau}$ be the transpose of the first column of $pr(A)$. We would like to start with the evidence of $b_i=a_i$. We have $$b_i=m-\sum^{l_1^{i-1}}_{s=i} a_{i-1,s}-\sum^n_{s=l_1^{i-1}} a_{i,s}-\sum^{i-1}_{s=2} pr(a_{s,i}),$$ $$a_i=m-\sum^n_{s=i}a_{i,s}-\sum^{i}_{s=2} pr(a_{s,i}),$$ since the sum over all entries of the $(i-1)$-th column of $B$ equals to $\sum^{l_1^{i-1}}_{s=i} a_{i-1,s}+\sum^n_{s=l_1^{i-1}} a_{i,s}$. However, this implies $$b_i-a_i=pr(a_{i,i})+\sum^{l_1^{i-1}-1}_{s=i}a_{i,s}-\sum^{l_1^{i-1}}_{s=i}a_{i-1,s}=pr(a_{i,i})-\epsilon(\mathbf a_{i-1},\mathbf a_{i})=0.$$ If $r>i$ we have $$b_r=\sum^{i-2}_{s=1} a_{s,r-1}-\sum^{i-1}_{s=2} pr(a_{s,r})+\begin{cases} a_{i-1,r-1}+a_{i,r-1}-\epsilon((\mathbf a_{i-1})^{\geq r},(\mathbf a_{i})^{\geq r}) &  \text{}\\ 
a_{i-1,r-1} & \text{}
\end{cases},$$
$$a_r=\sum^{i}_{s=1} a_{s,r-1}-\sum^{i}_{s=2} pr(a_{s,r})=\sum^{i}_{s=1} a_{s,r-1}-\sum^{i-1}_{s=2} pr(a_{s,r})-\begin{cases} \epsilon((\mathbf a_{i-1})^{\geq r},(\mathbf a_{i})^{\geq r}) &  \text{}\\ 
a_{i,r-1} & \text{}
\end{cases},$$ and thus the difference is once more zero. So by induction we can assume that the first row of $pr(B)$ consists of non-negative integers, and hence by our claim the first row of $pr(A)$ as well. Furthermore the sum over the entries in the last row of $B$  coincides with the sum over the entries in the last row of $A$ and thus $$m-\sum^n_{r=i}pr(a_{1,r})=\sum^i_{r=1}a_{r,n}.$$
\endproof
\end{prop}
At this point we take Lemma~\ref{index} and Remark~\ref{imrem} (2) up and emphasize that the action of $\tilde{f}_i$ on $pr(A)$, $A\in B^{m,i}$, is on one condition in the sense of the following lemma independent from the fact where $pr(A)$ lives. We will need this result several times in the remaining proofs.

\begin{lem}\label{guck}
Suppose that $A\in B^{m,i}$, $pr(A)\in B^{s,i}$ for some $s\geq m$ and  $\varphi_{i-1}(A)\neq 0$, then $_m\tilde{f}_i$ acts on $pr(A)$ and therefore $_s\tilde{f}_ipr(A)=\ _m\tilde{f}_ipr(A)$.
\proof
The sum over all entries in the $i$-th column of $pr(A)$ equals to the sum over all entries in the $(i-1)$-th column of $A$ and that is why \begin{flalign*}\sum^{i-1}_{r=1}pr(a_{1,r})+\sum^{n}_{r=i}pr(a_{i,r})&=m-\langle h_{i-1}, \wt(A) \rangle-\epsilon_{i-1}(A)&\\&=m-\varphi_{i-1}(A)<m.\end{flalign*}
\endproof
\end{lem}

\subsection{Main proofs} This section is dedicated to the verification of the conditions (1)-(3). Initially we remark that the cases $j=i-1$ and $j=i$ in part (3) of Proposition~\ref{prch} will be considered separately, which is the aim of the next proposition:
\begin{prop}\label{vor}
For $j=i-1,i$ and all $A$ in $B^{m,i}$ we have $$pr(\tilde{e}_{j}A)=\tilde{e}_{j+1}pr(A).$$
\proof
We presume $j=i-1$ and $\tilde{e}_{i-1}A\neq 0$, since the condition $\tilde{e}_{i-1}A=0$ forces $\epsilon_{i-1}(A)=pr(a_{i,i})=\epsilon_i(pr(A))=0$. So let $l_1,\cdots,l_t$ be the integers from step (1) if we apply the algorithm to $A$ and let $l'_1,\cdots,l'_{t'}$ be the integers obtained from step (1) if we apply our algorithm to $\tilde{e}_{i-1}A$. If $l'_1=l_1$, then the algorithm gives us immediately $\tilde{e}_ipr(A)=pr(\tilde{e}_{i-1}A)$. So suppose that $l_1>l'_1$ and let $d$ maximal such that $l'_d<l_1$. Then, using the definition of $l_1$, we get on the one hand $l'_{d+1}=l_1$ and on the other hand $$a_{i-1,l'_d+1}+\cdots+a_{i-1,l_1}-1=a_{i,l'_d}+\cdots+a_{i,l_1-1},$$ which means that $pr(\tilde{e}_{i-1}A)$ would not change if we skip $l'_d$. By repeating these arguments we can get rid of all $l'_d$ such that $l'_d<l_1$ and consequently we can calculate $pr(\tilde{e}_{i-1}A)$ by using the sequence $l'_{d+1}=l_1<l'_{d+2}=l_2<\cdots< l'_{t'}=l_t$. To be more accurate we can conclude $\tilde{e}_ipr(A)=pr(\tilde{e}_{i-1}A)$.\vskip 1pt Now let $j=i$ and in that additional separated case we will prove the required equation by induction on $i$. For the initial step we assume that $i=2$ and investigate the first two rows of $pr(A)$ where these are  
of one of the two following forms:
\begin{center}
$$\ytableausetup{boxsize=2.8em}
\begin{ytableau}
\bullet &  \epsilon_1(A) \\
a_{1,2} & a_{2,2}\\
\end{ytableau} \mbox{\quad or\quad }\ytableausetup{boxsize=2.8em}
\begin{ytableau}
\bullet & a_{1,2} \\
 y & \scriptstyle pr(a_{2,3})\\
\end{ytableau}\quad \mbox{ with } y=a_{1,2}+a_{2,2}-pr(a_{2,3})>a_{1,2},$$ 
\end{center}
whereas the first case appears if and only if either $q_-^1(A)>i$ or $q_-^1(A)=i$ and $a_{2,2}=pr(a_{2,3})$. In that case, since $\epsilon_1(A)\geq a_{1,2}$, we have $q_+^3(A)=2$ which means among other things that $\epsilon_3(pr(A))=0$ if $\epsilon_2(A)=a_{2,2}=0$. Furthermore, if $\tilde{e}_2A\neq 0$, we actually have $q_-^1(\tilde{e}_2A)>i$ provided $q_-^1(A)>i$ and $q_-^1(\tilde{e}_2A)=q_-((\mathbf a_1)^{>i},(\mathbf a_2)^{>i})$ provided $q_-^1(A)=i$ and $a_{2,2}=pr(a_{2,3})$. Thus $pr(\tilde{e}_2A)$ arises from $pr(A)$ by replacing $\epsilon_1(A)$ by $\epsilon_1(A)+1$ and $a_{2,2}$ by $a_{2,2}-1$, which proves the claim in that case. Otherwise the second case appears and there, because of $a_{2,2}>pr(a_{2,3})$, we have $\tilde{e}_2A\neq 0$, $q_+^3(A)=1$ and $q_-^1(\tilde{e}_2A)=i$. So the algorthm provides us the first two rows of $pr(\tilde{e}_2A)$: 
$$\ytableausetup{boxsize=2.8em}
\begin{ytableau}
\bullet+1 & a_{1,2} \\
 y-1 & \scriptstyle pr(a_{2,3})\\
\end{ytableau},$$ where the remaining entries coincide.
As a consequence we get in both cases $\tilde{e}_3pr(A)=pr(\tilde{e}_2A)$ so that we can devote our attention to the induction step by observing the element $B$ from the proof of Proposition~\ref{wise}. Suppose first that $\tilde{e}_iA\neq 0$ and let $B_{e_i}$ be the element obtained by same construction out of $\tilde{e}_iA$. We remember that the connection between $pr(A)$ and $pr(B)$ was $pr(A)=pr(B)pr(\mathbf a_i).$ 
By induction we can conclude among other things $$\epsilon_i(pr(B))=\epsilon_{i-1}(B)=\begin{cases} 
a_{i-1,i}+a_{i,i}-\epsilon((\mathbf a_{i-1})^{>i},(\mathbf{a_i})^{>i}) & \text{if $q^{i-1}_-(A)=q_-^{i-1}(\tilde{e_i}A)=i$},\\
a_{i-1,i}, & \text{ else. }
\end{cases}$$ The first opportunity forces $pr(a_{i,i})=a_{i-1,i}$, $pr(a_{i,i+1})=\epsilon((\mathbf a_{i-1})^{>i},(\mathbf{a_i})^{>i})<a_{i,i}$ where the second one forces $pr(a_{i,i})=\epsilon_{i-1}(A)\geq a_{i-1,i}$, $pr(a_{i,i+1})=a_{i,i}$. Moreover we have $q_+^{i+1}(pr(A))\in\{i,q_+^{i}(pr(B))\}$ and further we claim the following:
\begin{equation}\label{hun}q_+^{i+1}(pr(A))=q_+^{i}(pr(B))\mbox{ if and only if } q_-^{i-1}(A)=q_-^{i-1}(\tilde{e_i}A)=i\end{equation}
\textit{Proof of (\ref{hun}):}
We start by supposing $q_-^{i-1}(A)=q_-^{i-1}(\tilde{e_i}A)=i$ and get $$\epsilon_{i+1}(pr(A))\geq\epsilon_{i}(pr(B))-pr(a_{i,i})+pr(a_{i,i+1})=a_{i,i}>pr(a_{i,i+1}),$$ which implies $q_+^{i+1}(pr(A))\neq i$. For the converse direction let $q_-^{i-1}(A)=q_-^{i-1}(\tilde{e_i}A)=i$ be incorrect, then another easy estimation 
$$\epsilon_{i}(pr(B))-pr(a_{i,i})+pr(a_{i,i+1})=a_{i-1,i}-\epsilon_{i-1}(A)+a_{i,i}\leq a_{i,i}$$ implies $q_+^{i+1}(pr(A))=i$ and thus (\ref{hun}).\vskip 1pt Assume firstly that $q_-^{i-1}(A)=q_-^{i-1}(\tilde{e_i}A)=i$ which implies $B_{e_i}=\tilde{e}_{i-1}B$ and that the $i$-th row of $pr(\tilde{e}_iA)$ coincides with the $i$-th row of $pr(A)$.
According to (\ref{hun}) the gluing process commutes with the Kashiwara action. To be more precise $$pr(\tilde{e}_iA)=pr(B_{e_i})pr(\mathbf a_i)=(\tilde{e}_ipr(B))pr(\mathbf a_i)=\tilde{e}_{i+1}(pr(B)pr(\mathbf a_i))=\tilde{e}_{i+1}pr(A).$$  Lastly we assume that $q_-^{i-1}(A)=q_-^{i-1}(\tilde{e_i}A)=i$ is not fulfilled which implies immediately that $pr(\tilde{e}_iA)$ arises from $pr(A)$ by replacing $pr(a_{i,i})$ by $pr(a_{i,i})+1$ and $pr(a_{i,i+1})$ by $pr(a_{i,i+1})-1$. After that we finished our proof for all $A$ with $\tilde{e}_iA\neq 0$, since (\ref{hun}) implies exactly $q_+^{i+1}(pr(A))=i$. Although the ideas of the proof of the remaining case $\tilde{e}_iA=a_{i,i}=0$ are similar we will give it nevertheless for completeness. By the reason of $0=a_{i,i}\geq \epsilon((\mathbf a_{i-1})^{>i},(\mathbf{a_i})^{>i})$ we must necessarily have $\epsilon_{i}(B)=a_{i-1,i},pr(a_{i,i})=\epsilon_{i-1}(A)\geq a_{i-1,i}$, $pr(a_{i,i+1})=0$ and together with $\epsilon_{i}(pr(B))-pr(a_{i,i})=a_{i-1,i}-\epsilon_{i-1}(A)\leq 0$ we find $\epsilon_{i+1}(pr(A))=0.$
\endproof
\end{prop}
Hereafter we consider the remaining nodes:
\begin{prop}\label{abc}
The map $pr$ described in the algorithm satisfies the condition $$pr\circ\tilde{e}_{j}=\tilde{e}_{j+1}\circ pr \mbox{ for all $j\in\{1,\cdots,n-1\}$}.$$
\proof
By Proposition~\ref{vor} it is sufficient to verify the above stated equation for all $j<i-1$ and $j>i$, where we start by assuming that $j<i-1$ and for simplicity we set $q:=q^{j}_-(A)$. The basic idea of the proof is to compare permanently $pr(A)$ with $pr(\tilde{e}_{j}A)$ and reduce all assertions to the following claim:\vskip 6pt
\textbf{Claim 1:}\begin{enumerate}[i)]
\item Let $\tilde{e}_{j}A\neq 0$, then there exists an integer $z$, such that $pr(\tilde{e}_{j}A)$ arises out of $pr(A)$, if we replace $pr(a_{j+1,z})$ by $pr(a_{j+1,z})-1$ and $pr(a_{j+2,z})$ by $pr(a_{j+2,z})+1$
\item $q^{j+1}_-(pr(A))=z$
\end{enumerate}
Note that the claim will give us the proposition for all $j<i-1$, such that $\tilde{e}_{j}A\neq 0$. We want to emphasize here that in the proof of the claim we will also prove the statement of Proposition~\ref{abc} if $\tilde{e}_{j}A\neq 0$ is not satisfied. 
\vskip 3pt
\textit{Proof of Claim 1:}
For simplicity we denote by $\mathbf a^t=\ytableausetup{boxsize=1.4em}
\begin{ytableau}
a_{i} & \dots
&  a_{n} \\
\end{ytableau}$ the transpose of the $j$-th column of $A$ and by $\mathbf b^t=\ytableausetup{boxsize=1.4em}
\begin{ytableau}
b_{i} & \dots
&  b_{n} \\
\end{ytableau}$ the transpose of the $(j+1)$-th column of $A$. Further we will denote by $\mathbf c^t=\ytableausetup{boxsize=1.4em}
\begin{ytableau}
c_{i} & \dots
&  c_{n} \\
\end{ytableau}$ the transpose of the new column (\ref{newcolumn}) which we obtain after applying $(i-j-2)$ steps of our algorithm to $A$. For instance, if $j=i-2$, then $\mathbf c$ is precisely the $i$-th column of $A$. With the aim to obtain the $(j+2)$-th column of $pr(A)$, we apply our algorithm to the columns $\mathbf b$ and $\mathbf c$ and since there is no confusion we omit all superfluous indices and denote by $l_1,\cdots,l_t$ the integers described in the algorithm and suppose that $l_r< q \leq l_{r+1}$, where we understand again $l_0=i-1$. Denote by $(\mathbf {c'})^t=\ytableausetup{boxsize=1.4em}
\begin{ytableau}
c'_{i} & \dots
&  c'_{n} \\
\end{ytableau}$ the transpose of the $(j+2)$-th column of $pr(A)$ and by $(\mathbf {b'})^t=\ytableausetup{boxsize=1.4em}
\begin{ytableau}
b'_{i} & \dots
&  b'_{n} \\
\end{ytableau}$ the transpose of the new obtained column after $(i-j-1)$ steps of our algorithm. Hence these columns, expressed in terms of the entries of $\mathbf b$ and $\mathbf c$, are of the following form:
$$(\mathbf {b'})^t=\ytableausetup{boxsize=2.2em}
\begin{ytableau}
 b_i&
 b_{i+1}& \dots  & 
 b_{l_s-1}&
 b'_{l_s}&
b_{l_s+1}&
\dots&
\scriptstyle b_{n}+c_n
\\\end{ytableau} \quad s=1,\cdots,t-1$$

$$(\mathbf {c'})^t=\ytableausetup{boxsize=2.2em}
\begin{ytableau}
 x_{l_0}&
 c_{i}& \dots  & 
 c_{l_s-1}&
 x_{l_s}&
 c_{l_s+1}&
\dots&
 c_{n-1}
\\\end{ytableau}\quad s=1,\cdots,t-1,$$
whereby $x_{l_s}=b_{l_s+1}+\cdots+b_{l_{s+1}}-c_{l_s+1}-\cdots-c_{l_{s+1}-1}$ and $b'_{l_s}=b_{l_s}+c_{l_s}-x_{l_s}$.\vskip 3pt
With the goal to prove Claim 1 we need several minor results listed in Claim 1.1. \vskip 6pt
\textbf{Claim 1.1.:}
Let $s:=m_1=q_-(\mathbf a,\mathbf b'),\cdots,m_p$ be the integers obtained from the algorithm if we compare $\mathbf a$ with $\mathbf b'$ and let as before $l_r<q\leq l_{r+1}$. We presume further that $k$ is either \begin{equation}\label{coos}\min\{1\leq k\leq r|\sum^{l_{r+1}}_{r=l_k+1} b_r=\sum^{l_{r+1}-1}_{r=l_k} c_r\}\end{equation} or if the minium does not exist we set $k=r+1$. Then we have
\begin{enumerate}[i)]
\item $s\leq q$ and if $s<q$ then in fact $s\leq l_{k-1}$
\item $q\in\{m_1,\cdots,m_p\}$ and $\sharp\{x|l_{k-1}<m_x<q\}=0$
\item $\epsilon_{j+1}(pr(A))=\epsilon_j(A)$
\item $q^{j+1}_-(pr(A))\leq l_{k-1}+1$
\end{enumerate}
\textit{Proof of Claim 1.1.:}
Suppose $s\leq q$ in Claim 1.1.(i) is not fulfilled and let $l_p< s \leq l_{p+1}$. Then by observing the entries of $\mathbf b'$ we see that the sum $D:=a_i+\cdots+a_s+b'_s+\cdots+b'_n$ is of the following form:
$$D=\sum^s_{r=i}a_{r}+ \sum^{l_{p+1}}_{r=s}b_{r}+\sum^n_{r=l_{p+1}}c_{r}.$$But $$D\leq \sum^q_{r=i}a_{r}+ \sum^{l_{p+1}}_{r=q}b_{r}+\sum^n_{r=l_{p+1}}c_{r}\leq \sum^q_{r=i}a_{r}+ \sum^{l_{r+1}}_{r=q}b_{r}+\sum^n_{r=l_{r+1}}c_{r}=\sum^q_{r=i}a_r+\sum^n_{r=q}b'_r,$$ where the second last inequality is a consequence of the definition of $q$, particularly $\sum^{s-1}_{r=q}b_r\geq \sum^{s}_{r=q+1}a_r$ and the last inequality is by the definition of $l_{r+1}$, namely $\sum^{l_{p+1}-1}_{p=l_{r+1}}c_p\geq\sum^{l_{p+1}}_{p=l_{r+1}+1}b_p$. Consequently we obtain a contradiction to the definition of $s$ and thus \begin{equation}\label{qgrs} q\geq s=q_-(\mathbf a,\mathbf b').\end{equation} 
Before we start with the proof of the second statement in (i) we would like to emphasize the following result:
for all $j\in\{k,\cdots,r\}$ \begin{equation}\label{versch}\sum^{l_{r+1}}_{r=l_j+1} b_r=\sum^{l_{r+1}-1}_{r=l_j} c_r \mbox{ and }\ b'_{l_j}=b_{l_j} \mbox{ hold. }\end{equation}
Let us start by proving the first part of (\ref{versch}) by induction, where the initial step is by the choice of $k$ obvious. So assume that the first part of (\ref{versch}) holds for $j$. Using the definition of $l_j$ we must have $\sum^{s}_{r=l_j+1} b_r-\sum^{s-1}_{r=l_j+1} c_r\leq c_{l_j}$ for all $s>l_j$ and since $l_{j+1}$ is the ``place" where $\sum^{s}_{r=l_j+1} b_r-\sum^{s-1}_{r=l_j+1} c_r$ is maximal we have, together with the induction hypothesis, only one opportunity, namely $\sum^{l_{j+1}}_{r=l_{j}+1} b_r=\sum^{l_{j+1}-1}_{r=l_{j}} c_r.$ Hence the first part of (\ref{versch}) is proven. We proved also implicitly the second part which we can see as well as a corollary of the first part, namely we get for all $j\in\{k,\cdots,r\}$ \begin{flalign*}\sum^{l_{r+1}}_{r=l_{j}+1} b_r-\sum^{l_{r+1}-1}_{r=l_{j}} c_r=\sum^{l_{r+1}}_{r=l_{j+1}+1} b_r-\sum^{l_{r+1}-1}_{r=l_{j+1}} c_r=0\Rightarrow \sum^{l_{j+1}}_{r=l_{j}+1} b_r=\sum^{l_{j+1}-1}_{r=l_{j}} c_r\end{flalign*} which forces $b'_{l_{j}}=b_{l_{j}}$.
Because of (\ref{versch}) we verified part (i) of Claim 1.1. since the assumption $l_{k-1}<s<q$ would end in a contradiction, namely in $$D< D+\sum^q_{r=s+1}a_{r}-\sum^{q-1}_{r=s}b_{r}=D+\sum^q_{r=s+1}a_{r}-\sum^{q-1}_{r=s}b'_{r}=\sum^q_{r=i}a_r+\sum^n_{r=q}b'_r.$$
As a corollary of Claim 1.1. (i) we obtain that \begin{equation}\label{qqq}q\in\{m_1,\cdots,m_p\} \mbox{ and }\sharp\{x|l_{k-1}<m_x<q\}=0,\end{equation} because if $q=s$ we are done and if not we get with the definition of $m_2$ and similar calculations as in the proof of Claim 1.1. (i) that $q\geq m_2=q_-((\mathbf a)^{>s},(\mathbf b')^{>s})$ and in the case of $q>m_2$ we have $m_2\leq l_{k-1}$. If $q=m_2$ we are done and if not we repeat these arguments until we get (\ref{qqq}). 
For the completion of Claim 1.1. it remains to verify (iii) and (iv), whereas we start with \begin{equation}\label{zwe} \epsilon_{j+1}(pr(A))=\epsilon_{j}(A).\end{equation}
Note that the transpose of the $(j+1)$-th column of $pr(A)$ is given by 
$$\ytableausetup{boxsize=2.5em}
\begin{ytableau}
 z_{m_0}&
 b'_{i}& \dots  & 
 b'_{m_s-1}&
 z_{m_s}&
 b'_{m_s+1}&
\dots&
 b'_{n-1}
\\\end{ytableau}\quad s=1,\cdots,p-1,$$
whereby $z_{m_s}=a_{m_s+1}+\cdots+a_{m_{s+1}}-b'_{m_s+1}-\cdots-b'_{m_{s+1}-1}.$
Hence any sum $$\sum^h_{r=i} pr(a_{j+1,r})+\sum^n_{r=h} pr(a_{j+2,r}),$$ for some $i\leq h \leq n$ such that $m_{j-1}<h\leq m_j$ and $l_{p}<h\leq l_{p+1}$, is of the form \begin{equation}\label{asd}\sum^{m_j}_{r=i}a_r-\sum^{m_j-1}_{r=h}b'_{r}+\sum^{l_{p+1}}_{r=h}c'_r-\sum^{n}_{r=l_{p+1}+1}b_{r}\end{equation} and the expression $$\sum^h_{r=i} pr(a_{j+1,r})-\sum^{h-1}_{r=i} pr(a_{j+2,r})$$ can be written as
\begin{equation}\label{asd1}\sum^{m_j}_{r=i}a_r-\sum^{m_j-1}_{r=h}b'_{r}-\sum^{h-1}_{r=l_{p}+1}c'_r-\sum^{l_{p}}_{r=i}b_{r}.\end{equation} Assume that $h$ is minimal such that (\ref{asd}) is maximal, then since $b'_{k}\geq b_k$ for all $k=i,\cdots,n$ we have \begin{flalign*}\epsilon_{j+1}(pr(A))&=\sum^{m_j}_{r=i}a_r-\sum^{m_j-1}_{r=h}b'_{r}-\sum^{h-1}_{r=l_{p}+1}c'_r-\sum^{l_{p}}_{r=i}b_{r}&\\&\leq\sum^{m_j}_{r=i}a_r-\sum^{m_j-1}_{r=h}b_{r}+\begin{cases}-b_{i}\cdots-b_{l_p},&  \text{if $h=l_{p}+1$}\\
c_{l_{p+1}-1}+\cdots+c_{h-1}-b_{i}\cdots-b_{l_{p+1}},& \text{else}\end{cases}&\\&\leq \sum^{m_j}_{r=i}a_r-\sum^{m_j-1}_{r=i}b_{r}\leq \epsilon_j(A).\end{flalign*}The second last inequality is by the reason of $c_{h-1}+\cdots+c_{l_{p+1}-1}<b_{h}+\cdots+b_{l_{p+1}}$, which is valid by the definition of $l_{p+1}$ and $h-1\neq l_{p}$.
For the converse direction we investigate (\ref{asd1}) with $h=l_{k-1}+1$, whereby we can presume with (\ref{qqq}) that $q\in\{m_1,\cdots,m_p\}$, say $q=m_j$, and $l_{k-1} < h \leq l_{k}$, $m_{j-1} < h \leq m_j$. In addition we recall from the definition and (\ref{versch}) that $b'_{k}=b_k$ for $k=h,\cdots,q-1$ and obtain the reverse estimation, which will finish the proof of Claim 1.1. (iii):
\begin{flalign*}\epsilon_{j+1}(pr(A))&\geq\sum^{q}_{r=i}a_r-\sum^{q-1}_{r=h}b'_{r}-\sum^{h-1}_{r=l_{k-1}+1}c'_r-\sum^{l_{k-1}}_{r=i}b_{r}&\\&=\sum^{q}_{r=i}a_r-\sum^{q-1}_{r=h}b_{r}-b_{i}\cdots-b_{h-1}&\\&= \sum^{q}_{r=i}a_r-\sum^{q-1}_{r=i}b_{r}= \epsilon_j(A).\end{flalign*}
With these calculations we get among other things also $q_-^{j+1}(pr(A))\leq l_{k-1}+1$, because $\epsilon_{j+1}(pr(A))=\sum^{l_{k-1}+1}_{r=i} pr(a_{j+1,r})-\sum^{l_{k-1}}_{r=i} pr(a_{j+2,r})$ and $q_-^{j+1}(pr(A))$ is minimal with this property. \vskip 3pt
Now we return to the goal to convince ourselves from Claim 1 and fix some notation for $\tilde{e}_jA$. Let $\mathbf e_j$$\mathbf a$ and $\mathbf e_j$$\mathbf b$ respectively be the $j$-th and $(j+1)$-th column respectively of $\tilde{e}_{j}A$. We denote by $(\mathbf e_j\mathbf c')^t=\ytableausetup{boxsize=1.8em}
\begin{ytableau}
e_j c'_i & \dots
&  e_j c'_n \\
\end{ytableau}$ the transpose of the $(j+2)$-th column of $pr(\tilde{e}_{j}A)$ and by $(\mathbf e_j\mathbf b')^t=\ytableausetup{boxsize=1.8em}
\begin{ytableau}
e_j b'_i & \dots
&  e_j b'_n \\
\end{ytableau}$ we will denote the transpose of the new obtained column after applying $(i-j-1)$ steps of the algorithm to $\tilde{e}_{j}A$. For the purpose of determining $\mathbf e_j$$\mathbf c'$ we compare $\mathbf e_j$$\mathbf b$ with $\mathbf c$ and let $n_1,\cdots,n_y$ be the integers defined in the algorithm. Suppose again that $k$ is as in (\ref{coos}). If the minimum does not exist, then the integers do not change, i.e. $t=y$, $n_h=l_h$ for $h=1,\cdots,t$ and otherwise $k$ is the minimal integer such that $n_k\neq l_k$. A short calculation by using (\ref{versch}) shows that the new sequence of integers is given by $n_1=l_1,\cdots, n_{k-1}=l_{k-1},n_{k}=l_{r+1},\cdots, n_{k+t-r-1}=l_{t}$ and hence $$e_j c'_h=c'_h  \mbox{ for $h\neq l_{k-1}+1$ and } e_j c'_{l_{k-1}+1}=c'_{l_{k-1}+1}+1$$ and  $$e_j b'_h=b'_h \mbox{ for $h\neq l_{k-1},q$ and } e_j b'_{l_{k-1}}=b'_{l_{k-1}}-1,\ e_j b'_{q}=b'_{q}+1.$$
As a next step we compare the columns $\mathbf e_j\mathbf a$ with $\mathbf e_j\mathbf b'$ and determine the sequence of integers from step (1) of the algorithm, say $\overline{m}_1,\cdots,\overline{m}_x$. One can observe similar to (\ref{qgrs}) that $\overline{m}_1=q_-(\mathbf e_j\mathbf a,\mathbf e_j\mathbf b')\leq q$ and as a corollary we obtain again $q\in\{\overline{m}_1,\cdots,\overline{m}_x\}$. Accordingly we have the following situation:
$$\overline{m}_1<\overline{m}_2<\cdots<\overline{m}_{x_0}\leq l_{k-1} < \overline{m}_{x_1}<\cdots<\overline{m}_{x_n}=q<\cdots<\overline{m}_{x}$$
and since by (\ref{qqq}) there is no $1\leq h\leq p$ such that $l_{k-1}<m_h<q$ we have
$$m_1<m_2<\cdots<m_{j-1}\leq l_{k-1} < q=m_{j}<\cdots<m_p.$$
The integer $\overline{m}_{x_{n-1}}$ in the aforementioned sequence $\overline{m}_1<\cdots<\overline{m}_x$ has firstly the property $l_{k-1}<\overline{m}_{x_{n-1}}<q$ and secondly $\overline{m}_{x_{n-1}}$ is maximal with this property. Using the definition of $q$ we obtain $$b_{\overline{m}_{x_{n-1}}}+\cdots+b_{q-1}=a_{\overline{m}_{x_{n-1}}+1}+\cdots+a_q-1$$ and as a consequence we get with (\ref{versch}) that the resulting new obtained column (\ref{newcolumn}) does not change if we skip $\overline{m}_{x_{n-1}}$. Repeating these arguments we can get rid of all integers greater than $l_{k-1} $ and less than $q$ appearing in the sequence. Now it is obvious to see that we can replace the integer sequence $\overline{m}_1,\cdots,\overline{m}_x$ by $m_1,\cdots,m_p$ if we apply our algorithm to $\mathbf e_j\mathbf a$ and $\mathbf e_j\mathbf b'$. Thus we get two facts: the first fact is that the new obtained column at which we arrive by applying the algorithm to $\mathbf e_j$$\mathbf a$ and $\mathbf e_j$$\mathbf b'$ is the same as the one if we apply the algorithm to $\mathbf a$ and $\mathbf b'$. The second fact is that the $(j+1)$-th column of $pr(\tilde{e}_{j}A)$ is almost the same as the $(j+1)$-th column of $pr(A)$ except the $(l_{k-1}+1)$-th entry is one smaller. Consequently, we proved part (i) of Claim 1.
Now part (ii) of Claim 1 is also proven since (\ref{zwe}) forces $q^{j+1}_-(pr(A))\geq z=l_{k-1}+1$, because otherwise we get $\epsilon_{j+1}(pr(\tilde{e}_{j}A))=\epsilon_{j+1}(pr(A))$ which is a contradiction to $$\epsilon_{j+1}(pr(\tilde{e}_{j}A))=\epsilon_{j}(\tilde{e}_{j}A)=\epsilon_{j}(A)-1=\epsilon_{j+1}(pr(A))-1.$$
Consequently we proved our proposition for all $j<i-1$. From now on we would like to show \begin{equation}\label{dri}\tilde{e}_{j+1}pr(A)=pr(\tilde{e}_{j}A)\mbox{ for all $j>i$. }\end{equation}If $A$ is any element in $B^{m,i}$ such that $\varphi_1(A)=\cdots=\varphi_{i-1}(A)=0$, then it is easy to see that the image under $pr$ is given by
$$pr(a_{r,s})=\begin{cases} 
\epsilon_{r-1}(A), & \text{if $s=i, r>1$}\\
m-\sum^n_{p=i} a_{1,p}, & \text{if $s=i, r=1$}\\
a_{r,s-1},& \text{ else}.
\end{cases}$$
Note that $\tilde{e}_jA=0$ implies $\tilde{e}_{j+1}pr(A)=0$ and otherwise, using the Stembridge axioms which are fulfilled by Theorem~\ref{nakmon} for all elements in $B^{m,i}$, we obtain
$\varphi_1(\tilde{e}_jA)=\cdots=\varphi_{i-1}(\tilde{e}_jA)=0$ and therefore by applying the algorithm to $\tilde{e}_jA$ and comparing $pr(\tilde{e}_{j}A)$ with $pr(A)$ we have $\tilde{e}_{j+1}pr(A)=pr(\tilde{e}_{j}A)$. So we proved our proposition for all $A$, such that $\varphi_1(A)=\cdots=\varphi_{i-1}(A)=0$. Now let $A$ be arbitrary and write $\wt(A)$ as a linear combination of simple roots $\wt(A) = \sum_{j \in I} k_j \alpha_j\in \sum_{j \in I} \Q \alpha_j $, and define
\[ \ts
  \height(\wt(A)) := \sum_{j \in I} k_j.
\]
Our proof will proceed by induction on $\lceil \height(\wt(A)) \rceil$, whereby $\lceil \cdot \rceil$ denotes the ceiling function. If the height is minimal we have the lowest weight element in $B^{m,i}$ which satisfies obiously (\ref{dri}). If $\varphi_1(A)=\cdots=\varphi_{i-1}(A)=0$, then we are done by the above considerations and if not let $1\leq l \leq i-1$ be any integer so that $\varphi_l(A)\neq0$. By induction we gain $\tilde{e}_{j+1}pr(\tilde{f}_lA)=pr(\tilde{e}_{j}\tilde{f}_lA)$ and by earlier calculations, together with Lemma~\ref{index}, Lemma~\ref{guck} and Proposition~\ref{vor}, we can verify $pr(A)=pr(\tilde{e}_l\tilde{f}_{l}A)=\tilde{e}_{l+1}pr(\tilde{f}_{l}A)\Rightarrow \tilde{f}_{l+1}pr(A)=pr(\tilde{f}_lA)$. Thus, by using the Stembridge axioms, we can finish the proof of $(\ref{dri})$ since
$$\epsilon_{j+1}(pr(A))=\epsilon_{j+1}(\tilde{f}_{l+1}pr(A))=\epsilon_{j+1}(pr(\tilde{f}_lA))=\epsilon_{j}(\tilde{f}_lA)=\epsilon_{j}(A)$$
and \begin{flalign*}\tilde{f}_{l+1}\tilde{e}_{j+1}pr(A)&=\tilde{e}_{j+1}\tilde{f}_{l+1}pr(A)=\tilde{e}_{j+1}pr(\tilde{f}_{l}A)=pr(\tilde{e}_{j}\tilde{f}_lA)&\\&=pr(\tilde{f}_{l}\tilde{e}_jA)=\tilde{f}_{l+1}pr(\tilde{e}_{j}A).\end{flalign*}
\endproof
\end{prop}
At this point we are in position to state our main theorem:
\begin{thm}\label{mainthm0}
The map $pr$ described in the algorithm is Schützenberger's promotion operator.
\proof
Let $A\in B^{m,i}$, then we erase all arrows with colour $n$ and denote by $Z_{(1,\cdots,n-1)}(A)$ the connected component containing $A$. Let $B$ be the $\{1,\cdots,n-1\}$ highest weight element. Then by an immediate inspection of the definiton we must have $b_{r,s}=0$ for all $(r,s)$ except $(r,s)=(i,n)$, which in particular means $$pr(b_{r,s})=\begin{cases} 
m-b_{i,n}, & \text{if $r=1, s=i$}\\
0,& \text{ else}.
\end{cases}$$
Thus we have $pr(\tilde{e}_{i_1}\cdots \tilde{e}_{i_s}A)=pr(B)\in B^{m,i}$, with some $i_1,\cdots i_s\in \{1,\cdots,n-1\}$. We claim actually that $pr(A)$ lives in $B^{m,i}$ and we will prove this statement by induction on $n:=\sharp\{i_r|i_r=i-1\}$. Suppose that $pr(A)\in B^{p,i}$ for some $p\geq m$. If $n=0$ we obtain by Lemma~\ref{index}
$$pr(A)=\ _p\tilde{f}_{i_s+1}\cdots _p\tilde{f}_{i_1+1}pr(B)=\ _m\tilde{f}_{i_s+1}\cdots _m\tilde{f}_{i_1+1}pr(B)\in B^{m,i},$$ which proves the initial step. Now we assume that $l=\min\{1\leq l \leq s|i_l=i-1\}$ and $pr(\tilde{e}_{i_{l+1}}\cdots \tilde{e}_{i_s}A)\in B^{p,i}$ for some $p\geq m$. Then we get $$pr(\tilde{e}_{i_{l+1}}\cdots \tilde{e}_{i_s}A)= \ _p\tilde{f}_{i_{l}+1}\cdots _p\tilde{f}_{i_1+1}pr(B)=\ _p\tilde{f}_{i_{l}+1}\ _m\tilde{f}_{i_{l-1}+1}\cdots _m\tilde{f}_{i_1+1}pr(B).$$ 
So by the induction hypothesis it is sufficient to prove that \begin{equation}\label{bum}\ _p\tilde{f}_{i_{l}+1}\ _m\tilde{f}_{i_{l-1}+1}\cdots _m\tilde{f}_{i_1+1}pr(B)\in B^{m,i},\end{equation}
but since $\varphi_{i-1}(\tilde{e}_{i_{l}}\cdots \tilde{e}_{i_s}A)\neq 0$ we can conclude with Lemma~\ref{guck} that (\ref{bum}) holds: $$_p\tilde{f}_{i_{l}+1}\ _m\tilde{f}_{i_{l-1}+1}\cdots _m\tilde{f}_{i_1+1}pr(B)=\ _m\tilde{f}_{i_{l}+1}\ _m\tilde{f}_{i_{l-1}+1}\cdots _m\tilde{f}_{i_1+1}pr(B)\in B^{m,i}.$$ 
%However,\begin{flalign*}\varphi_{i}(_m\tilde{f}_{i_{l-1}+1}\cdots _m\tilde{f}_{i_1+1}pr(B))&=\varphi_{i}(pr(\tilde{e}_{i_l}\cdots \tilde{e}_{i_s}A))&\\&=m-\langle h_{i-1}, \wt(\tilde{e}_{i_l}\cdots \tilde{e}_{i_s}A) \rangle-\epsilon_{i-1}(\tilde{e}_{i_l}\cdots \tilde{e}_{i_s}A)&\\&=m-\varphi_{i-1}(\tilde{e}_{i_l}\cdots \tilde{e}_{i_s}A)&\\&=m-\varphi_{i-1}(\tilde{e}_{i_{l+1}}\cdots \tilde{e}_{i_s}A))-1<m.\end{flalign*}
According to that we have the well-definedness of $pr$, i.e. $pr:B^{m,i}\longrightarrow B^{m,i}$. The condition (1) of Proposition~\ref{prch} is obviously fulfilled by construction and Proposition~\ref{wise} and condition (3) is exactly Proposition~\ref{abc} and the following simple calculation: 
by the reason of condition (1) of Proposition~\ref{prch} and part (1) of Definition~\ref{abcrystal} we can assume without loss of generality that $\tilde{f}_jA\neq 0$ and thus $$pr(A)=pr(\tilde{e}_j\tilde{f}_jA)=\tilde{e}_{j+1}pr(\tilde{f}_jA).$$
So the proof of part (2) of Proposition~\ref{prch} will finish our main theorem.
Note that for the bijectivity it is enough to prove the surjectivity. So let $A\in B^{m,i}$ be an arbitrary element and let $B$ the highest weight element in $Z_{(2,\cdots,n)}(A)$. Then it is obvious to see that $B$ has the property $b_{r,s}=0$ if $(r,s)\neq(1,i)$ and according to this $B$ has a pre-image, say $C$. For instance one can choose $C$ as follows: 
$$c_{r,s}=\begin{cases} 0,& \text{$(r,s)\neq(i,n)$}\\
m-b_{1,i}, & \text{if $r=i, s=n$}.
\end{cases}$$ Therefore, since $B=\tilde{e}_{i_1}\cdots \tilde{e}_{i_s}A$ with $i_1\cdots i_s\in \{2,\cdots,n\}$, we have $$A=\tilde{f}_{i_s}\cdots \tilde{f}_{i_1}B=\tilde{f}_{i_s}\cdots \tilde{f}_{i_1}pr(C)=pr(\tilde{f}_{i_s-1}\cdots \tilde{f}_{i_1-1}C).$$
\endproof
\end{thm}
\begin{rem}
If we follow the results from \cite{BST10} we can compute the inverse map of $pr$ by composing $n$ times $pr$. In particular $$pr^{-1}=pr^n.$$
\end{rem}
We would like to finish our paper with drawing a KR-crystal graph of type $A^{(1)}_2$.
\newpage
\begin{ex}
The KR-crystal $B^{3,2}$ of type $A^{(1)}_2$ looks as follows:
\end{ex}

\begin{center}
\begin{tikzpicture}[>=latex,line join=bevel,]
\node (1+1+1+2+3+3) at (140bp,360bp) [draw,draw=none] {${\def\lr#1{\multicolumn{0}{|@{\hspace{.6ex}}c@{\hspace{.6ex}}|}{\raisebox{-.3ex}{$#1$}}}\raisebox{-.6ex}{$\begin{array}[b]{ccc}\cline{1-1}\cline{2-2}\lr{0}&\lr{2}\\\cline{1-1}\cline{2-2}\end{array}$}}$};
  \node (1+1+2+2+3+3) at (80bp,274bp) [draw,draw=none] {${\def\lr#1{\multicolumn{1}{|@{\hspace{.6ex}}c@{\hspace{.6ex}}|}{\raisebox{-.3ex}{$#1$}}}\raisebox{-.6ex}{$\begin{array}[b]{ccc}\cline{1-1}\cline{2-2}\lr{1}&\lr{1}\\\cline{1-1}\cline{2-2}\end{array}$}}$};
  \node (1+2+2+2+3+3) at (19bp,188bp) [draw,draw=none] {${\def\lr#1{\multicolumn{1}{|@{\hspace{.6ex}}c@{\hspace{.6ex}}|}{\raisebox{-.3ex}{$#1$}}}\raisebox{-.6ex}{$\begin{array}[b]{ccc}\cline{1-1}\cline{2-2}\lr{2}&\lr{0}\\\cline{1-1}\cline{2-2}\cline{3-3}\end{array}$}}$};
  \node (1+1+1+3+3+3) at (201bp,274bp) [draw,draw=none] {${\def\lr#1{\multicolumn{1}{|@{\hspace{.6ex}}c@{\hspace{.6ex}}|}{\raisebox{-.3ex}{$#1$}}}\raisebox{-.6ex}{$\begin{array}[b]{ccc}\cline{1-1}\cline{2-2}\lr{0}&\lr{3}\\\cline{1-1}\cline{2-2}\cline{3-3}\end{array}$}}$};
  \node (2+2+2+3+3+3) at (37bp,16bp) [draw,draw=none] {${\def\lr#1{\multicolumn{1}{|@{\hspace{.6ex}}c@{\hspace{.6ex}}|}{\raisebox{-.3ex}{$#1$}}}\raisebox{-.6ex}{$\begin{array}[b]{ccc}\cline{1-1}\cline{2-2}\lr{3}&\lr{0}\\\cline{1-1}\cline{2-2}\cline{3-3}\end{array}$}}$};
  \node (1+2+2+3+3+3) at (74bp,102bp) [draw,draw=none] {${\def\lr#1{\multicolumn{1}{|@{\hspace{.6ex}}c@{\hspace{.6ex}}|}{\raisebox{-.3ex}{$#1$}}}\raisebox{-.6ex}{$\begin{array}[b]{ccc}\cline{1-1}\cline{2-2}\lr{2}&\lr{1}\\\cline{1-1}\cline{2-2}\cline{3-3}\end{array}$}}$};
  \node (1+1+2+2+2+3) at (19bp,360bp) [draw,draw=none] {${\def\lr#1{\multicolumn{1}{|@{\hspace{.6ex}}c@{\hspace{.6ex}}|}{\raisebox{-.3ex}{$#1$}}}\raisebox{-.6ex}{$\begin{array}[b]{ccc}\cline{1-1}\cline{2-2}\lr{1}&\lr{0}\\\cline{1-1}\cline{2-2}\cline{3-3}\end{array}$}}$};
  \node (1+1+1+2+2+2) at (36bp,532bp) [draw,draw=none] {${\def\lr#1{\multicolumn{1}{|@{\hspace{.6ex}}c@{\hspace{.6ex}}|}{\raisebox{-.3ex}{$#1$}}}\raisebox{-.6ex}{$\begin{array}[b]{ccc}\cline{1-1}\cline{2-2}\lr{0}&\lr{0}\\\cline{1-1}\cline{2-2}\cline{3-3}\end{array}$}}$};
  \node (1+1+1+2+2+3) at (74bp,446bp) [draw,draw=none] {${\def\lr#1{\multicolumn{1}{|@{\hspace{.6ex}}c@{\hspace{.6ex}}|}{\raisebox{-.3ex}{$#1$}}}\raisebox{-.6ex}{$\begin{array}[b]{ccc}\cline{1-1}\cline{2-2}\lr{0}&\lr{1}\\\cline{1-1}\cline{2-2}\cline{3-3}\end{array}$}}$};
  \node (1+1+2+3+3+3) at (140bp,188bp) [draw,draw=none] {${\def\lr#1{\multicolumn{1}{|@{\hspace{.6ex}}c@{\hspace{.6ex}}|}{\raisebox{-.3ex}{$#1$}}}\raisebox{-.6ex}{$\begin{array}[b]{ccc}\cline{1-1}\cline{2-2}\lr{1}&\lr{2}\\\cline{1-1}\cline{2-2}\cline{3-3}\end{array}$}}$};
  
  \draw [red,->] (1+2+2+2+3+3) ..controls (32bp,162bp) and (40bp,148bp)  .. (48bp,136bp) .. controls (50bp,133bp) and (53bp,129bp)  .. (1+2+2+3+3+3);
  \pgfsetstrokecolor{black}
  \draw (57bp,145bp) node {$2$};
  \draw [blue,->] (1+1+1+2+2+3) ..controls (57bp,425bp) and (52bp,418bp)  .. (48bp,412bp) .. controls (42bp,403bp) and (36bp,394bp)  .. (1+1+2+2+2+3);
  \draw (57bp,403bp) node {$1$};
  \draw [blue,->] (1+1+1+2+3+3) ..controls (123bp,339bp) and (118bp,332bp)  .. (114bp,326bp) .. controls (108bp,317bp) and (101bp,307bp)  .. (1+1+2+2+3+3);
  \draw (123bp,317bp) node {$1$};
  \draw [black,<-] (1+1+1+2+2+3) ..controls (76bp,384bp) and (78bp,320bp)  .. (1+1+2+2+3+3);
  \draw (85bp,360bp) node {$0$};
  \draw [black,<-] (1+2+2+2+3+3) ..controls (23bp,141bp) and (25bp,111bp)  .. (28bp,86bp) .. controls (30bp,67bp) and (33bp,46bp)  .. (2+2+2+3+3+3);
  \draw (37bp,102bp) node {$0$};
  \draw [blue,->] (1+1+2+3+3+3) ..controls (118bp,159bp) and (104bp,140bp)  .. (1+2+2+3+3+3);
  \draw (123bp,145bp) node {$1$};
  \draw [blue,->] (1+1+2+2+3+3) ..controls (59bp,253bp) and (53bp,246bp)  .. (48bp,240bp) .. controls (42bp,231bp) and (36bp,222bp)  .. (1+2+2+2+3+3);
  \draw (57bp,231bp) node {$1$};
  \draw [black,<-] (1+1+2+2+2+3) ..controls (19bp,298bp) and (19bp,234bp)  .. (1+2+2+2+3+3);
  \draw (28bp,274bp) node {$0$};
  \draw [black,<-] (1+1+1+2+2+2) ..controls (31bp,493bp) and (30bp,477bp)  .. (28bp,462bp) .. controls (25bp,431bp) and (22bp,396bp)  .. (1+1+2+2+2+3);
  \draw (37bp,446bp) node {$0$};
  \draw [red,->] (1+1+1+2+3+3) ..controls (160bp,331bp) and (173bp,313bp)  .. (1+1+1+3+3+3);
  \draw (184bp,317bp) node {$2$};
  \draw [red,->] (1+1+1+2+2+2) ..controls (49bp,503bp) and (56bp,485bp)  .. (1+1+1+2+2+3);
  \draw (66bp,489bp) node {$2$};
  \draw [red,->] (1+1+2+2+3+3) ..controls (96bp,248bp) and (106bp,234bp)  .. (114bp,222bp) .. controls (116bp,219bp) and (119bp,215bp)  .. (1+1+2+3+3+3);
  \draw (123bp,231bp) node {$2$};
  \draw [black,<-] (1+1+2+2+3+3) ..controls (78bp,212bp) and (76bp,148bp)  .. (1+2+2+3+3+3);
  \draw (85bp,188bp) node {$0$};
  \draw [red,->] (1+1+2+2+2+3) ..controls (32bp,334bp) and (39bp,319bp)  .. (48bp,308bp) .. controls (51bp,305bp) and (53bp,301bp)  .. (1+1+2+2+3+3);
  \draw (57bp,317bp) node {$2$};
  \draw [blue,->] (1+2+2+3+3+3) ..controls (62bp,73bp) and (55bp,55bp)  .. (2+2+2+3+3+3);
  \draw (68bp,59bp) node {$1$};
  \draw [blue,->] (1+1+1+3+3+3) ..controls (181bp,245bp) and (168bp,227bp)  .. (1+1+2+3+3+3);
  \draw (186bp,231bp) node {$1$};
  \draw [red,->] (1+1+1+2+2+3) ..controls (96bp,417bp) and (110bp,398bp)  .. (1+1+1+2+3+3);
  \draw (120bp,403bp) node {$2$};
  \draw [black,<-] (1+1+1+2+3+3) ..controls (140bp,298bp) and (140bp,234bp)  .. (1+1+2+3+3+3);
  \draw (149bp,274bp) node {$0$};
\end{tikzpicture}
\end{center}
% References
%%%%%%%%%%%%%%%%%%%%%%%%%%%%%%%%%%%%%%%%%%%%%%%%%%%%%%%%%%%%%%%%%%%
\bibliographystyle{alpha}
\bibliography{crystalstructureonpolytop-biblist}
\end{document}